\documentclass[preprint,11pt]{elsarticle}
\usepackage{geometry}
\usepackage{setspace}
\usepackage[usestackEOL]{stackengine}
\usepackage{rotating}
\usepackage{lipsum}
\usepackage[utf8]{inputenc}
\usepackage[T1]{fontenc}
\usepackage{float}
\usepackage{textcomp}
\usepackage{pgf}
\usepackage{tikz}
\usetikzlibrary{arrows,automata}
\usepackage{amsmath}

\usepackage{multirow}
\usepackage{algorithm}
\usepackage[noend]{algpseudocode}
\usepackage{subcaption}
\usepackage{verbatim}
\usepackage{changepage}
\usepackage{mathtools}
\usepackage{amsfonts}
\usepackage{amsthm}
\usepackage{empheq}
\usepackage[shortlabels]{enumitem}
\usepackage{amsmath, amssymb}
\usepackage{tabularray} 
\newcommand{\edit}[1]{{\color{black} #1}}
\newcommand{\revise}[1]{{\color{black} #1}}

\renewcommand{\eqref}[1]{\textup{{\normalfont(\ref{#1}}\normalfont)}}

\biboptions{authoryear,longnamesfirst}

\usepackage{tikz}
\usepackage{everypage}
\usepackage{hyperref}
\usepackage{xcolor}
\AddEverypageHook{%
  \begin{tikzpicture}[remember picture, overlay]
    \node[anchor=west, text=gray, font=\small] 
      at ([xshift=1in, yshift=-2cm]current page.north west)
      {\href{https://doi.org/10.1016/j.ejor.2025.06.025}{https://doi.org/10.1016/j.ejor.2025.06.025 [Paper accepted in EJOR/Final version]}};
  \end{tikzpicture}
}

\usepackage[finalnew]{trackchanges} 
\addeditor{NG}
\newtheorem{prop}{Proposition}
\newtheorem{cor}{Corollary}
\newtheorem{assum}{Assumption}

\newtheorem{rem}{Remark}
\newtheorem{eg}{Example}
\newtheorem{lemma}{Lemma}
 \usepackage{todonotes,bm}
\allowdisplaybreaks
\journal{European Journal of Operational Research}

\begin{document}

\begin{frontmatter}

\title{Data-driven Joint Optimization of Maintenance and Spare Parts Provisioning: A Distributionally Robust  Approach}
\author[label2]{Heraldo Rozas $^{*}$}
\author[label1]{Weijun Xie}
\author[label1]{Nagi Gebraeel}
\author[label3]{Stephen Robinson}

\affiliation[label2]{organization={Electrical Engineering Department, University of Chile and Instituto Sistemas Complejos de Ingeniería (ISCI)},
            addressline={2007 Tupper}, 
            city={Santiago},
            postcode={8330111},  
            country={Chile}}
\affiliation[label1]{organization={H. Milton Stewart School of Industrial and Systems Engineering, Georgia Institute of Technology},
            addressline={765 Ferst Drive}, 
            city={Atlanta},
            postcode={30332}, 
            state={GA},
            country={USA}} 
\affiliation[label3]{organization={Department of Mechanical and Aerospace Engineering, University of California Davis}, 
            city={Davis},
            postcode={95616}, 
            state={CA},
            country={USA}} 

\begin{abstract}
This paper investigates the joint optimization of condition-based maintenance and spare provisioning, incorporating insights obtained from sensor data. Prognostic models estimate components' remaining lifetime distributions (RLDs), which are integrated into an optimization model to coordinate maintenance and spare provisioning. The existing literature addressing this problem assumes that prognostic models provide accurate estimates of RLDs, thereby allowing a direct adoption of Stochastic Programming or Markov Decision Process methodologies. Nevertheless, this assumption often does not hold in practice since the estimated distributions can be inaccurate due to noisy sensors or scarcity of training data. To tackle this issue, we develop a Distributionally Robust Chance Constrained (DRCC) formulation considering general discrepancy-based ambiguity sets that capture potential distribution perturbations of the estimated RLDs. The proposed formulation admits a Mixed-Integer Linear Programming (MILP) reformulation, where explicit formulas are provided to simplify the general discrepancy-based ambiguity sets. Finally, for the numerical illustration, we test a type-$\infty$ Wasserstein ambiguity set and derive closed-form expressions for the parameters of the MILP reformulation. The efficacy of our methodology is showcased in a wind turbine case study, where the proposed DRCC formulation outperforms other benchmarks based on stochastic programming and robust optimization.
\end{abstract}



\begin{keyword}
Maintenance \sep Spare Parts Inventory \sep Distributionally Robust Chance Constrained Programs\sep Condition Monitoring 


\end{keyword}
\end{frontmatter} 
 ${}^{*}$ Corresponding author. 
 \textit{E-mail addresses:} heraldo.rozas@ug.uchile.cl (H.Rozas), wxie@gatech.edu (W.Xie),  nagi.gebraeel@isye.gatech.edu(N.Gebraeel), stephen.k.robinson@ucdavis.edu (S.Robinson)
 
\section{Introduction} \label{sec:intro}
\subsection{Motivation}

Many industrial enterprises are involved in various types of digital transformation initiatives. Some of these initiatives involve installing sensors and automation technologies that generate large volumes of data. The data in and of itself is not actionable unless appropriate data analytic models are applied to extract valuable insights. Data analytics improves productivity and profitability, reduces waste and environmental impact, and provides the required visibility to mitigate unexpected disruptions. Maintenance, Repairs, and Operations (MRO) represent an important area where data analytics is extensively applied within the industrial sector, driving significant improvements in efficiency, reliability, and cost savings \citep{olsen2020industry}. Sensors are used to monitor the condition of industrial equipment and machines. Condition Monitoring (CM) data from these sensors is processed using analytic models and algorithms to help companies determine more precisely when equipment maintenance and part replacement are necessary, i.e., condition-Based Maintenance (CBM). Advance knowledge of scheduled repairs and/or part replacements can also be used to manage spare parts inventory better, potentially saving companies millions of dollars in reduced downtime, inventory cost reduction, and increased throughput.

Prognostic modeling is a special class of predictive analytics solutions that uses CM data from sensors to predict the remaining useful life of partially degraded components. Some popular prognostic modeling frameworks, such as \mbox{\citep{gebraeel2005residual}}, have successfully managed to combine population-based reliability characteristics with component-specific CM data to compute and continuously update Remaining Life Distributions (RLDs) of machine components in real-time. These distributions contain critical information that can be leveraged to compute efficient CBM schedules. They can also be leveraged to efficiently plan just-in-time provisioning of spare parts \citep{wang2015prognostics,zheng2021joint}. 

Despite being a compelling proposition, there is still a lack of a robust framework that unifies the three components of this problem, i.e., prognostic modeling, CBM, and spare parts inventory management. The demand for spare parts is driven by repair operations, which are planned based on a CBM schedule. An accurate CBM schedule can be computed using RLDs from the prognostic models. Prognostic models leverage noisy component-specific CM data to compute RLDs. The effectiveness of this value chain greatly depends on the accuracy of the components' RLDs. It also depends on how well these noisy RLDs can be integrated within a robust optimization framework that jointly optimizes CBM and spare parts provisioning. This paper focuses on the joint optimization of maintenance schedules and spare parts provisioning, incorporating insights obtained from prognostic modeling. 

\subsection{Related work}\label{subsec:literature_review}

The optimization of maintenance schedules and spare inventory for industrial components has been traditionally studied separately. \edit{Recent review papers focusing separately on maintenance and spare provisioning have been published by \citep{arts2024fifty} and \citep{hu2018or}, respectively.} Maintenance models assume that spare parts are always available in stock, and so make decisions based solely on lifetime or remaining lifetime distributions \citep{barlow1960optimum,jardine1997optimal,drent2023real,yildirim2016sensor}. Spare inventory models, on the other hand, usually rely on time-series methodologies to predict the demand for spare parts based on the historical spare demand without considering the impact of maintenance decisions \citep{li2011bayesian,van2019forecasting,howard2015inventory}.  Nevertheless, optimizing these two problems in isolation or a sequential fashion can lead to sub-optimal solutions due to the tight coupling between spare inventory and maintenance decisions \citep{van2013joint}. This motivates the development of new methodologies to optimize maintenance schedules and spare provisioning jointly.

The joint optimization of maintenance and spare provisioning have been investigated under varied approaches. A thorough survey is presented in \citep{van2013joint}. The first papers on this topic focused on threshold-based policies using reliability distributions. For example, \cite{kabir1996stocking} developed a joint policy, known as the (s, S, t) stocking policy, which minimizes the long-run total cost of replacement and spare parts inventory. In this policy, the reorder point is denoted by $s$, the order-up-to level is denoted by $S$, and a specified age for preventive maintenance is denoted by $t$. The optimal policy was obtained through simulation. This work was extended in \citep{sarker2000optimization} by incorporating additional factors, such as different costs associated with different inventory levels, the number of failures, and the placement of emergency orders when the inventory level falls abruptly. \cite{hu2008joint} also extended this model by implementing a genetic algorithm to reduce the objective value and the running time. A notable limitation of this line of work is that the models rely on reliability information instead of CM data. This can lead to conservative policies that usually recommend frequent, unnecessary repairs to protect against unexpected failures \citep{yildirim2016sensor}.

The widespread adoption of Internet of Things (IoT) technologies has boosted the implementation of CBM strategies in industrial applications \citep{olsen2020industry}. As a result, recent papers have investigated using CM strategies to jointly optimize maintenance and spare inventory \citep{rausch2010joint,
zhang2017joint, keizer2017joint,feng2023managing}. \cite{rausch2010joint} analyzed a single-unit system with a spare part inventory pool and applied bi-objective optimization to minimize the spare part inventory and the expected total operating cost. \cite{zhang2017joint} and \cite{keizer2017joint} studied a system with multiple identical units by stochastic decision processes. They developed condition-based policies for maintenance and spare provisioning relying on genetic algorithms \citep{zhang2017joint} and Markov Decision Process (MDP) \citep{keizer2017joint}. \cite{feng2023managing} investigated the joint optimization of production, inventory, and maintenance for a group of identical machines by continuous-time MDP. These existing models \citep{zhang2017joint, keizer2017joint,rausch2010joint,feng2023managing} generated control policies that trigger maintenance activities when a degradation signal reaches a certain threshold. However, such decisions are based solely on the current conditions of the components without considering predictions of the evolution of components' degradation. Predicting components' remaining lifetime can be critical for planning future repair activities and thus plays a key role in the joint optimization problem.

Failure prognostic algorithms enable the estimation of components' RLDs based on CM data. The development of such algorithms has been prominent during the last decades--\cite{elattar2016prognostics} presents a survey on this topic. Estimated RLDs can provide advanced knowledge to jointly optimize maintenance activities and spare parts provisioning. However, the use of prognostics for jointly optimizing maintenance and spare inventory has not been studied extensively. The existing works on solving the joint optimization problem using prognostic results can be classified depending on the size of the system in the study, resulting in single-unit and multi-unit system approaches. Single-unit system approaches \citep{elwany2008sensor,wang2015prognostics,tang2024joint} deal with a system composed of a single component. Implementing these methodologies in industrial applications is limited as real-world systems often have multiple components. Multi-unit system approaches \citep{dieterman2019joint, zheng2021joint, shi2022stochastic,zheng2023joint,feng2023managing,zheng2024joint}, on the other hand, analyze systems with multiple components. However, the existing multi-unit approaches focus on systems with multiple \textit{identical} components. In this context, the joint optimization problem has been formulated and solved using different methodologies, including threshold-based policies \citep{dieterman2019joint}, MDP \citep{zheng2021joint,zheng2023joint,feng2023managing}, reinforcement learning \citep{zheng2024joint}, and stochastic programming \citep{shi2022stochastic}. The existing joint optimization models using prognostic results for multi-unit systems have two main limitations: 1) They assume that the RLDs estimated by the prognostic algorithms \textit{accurately} characterize the uncertainty of the components' failure times. This assumption allows for a direct integration of the estimated RLDs into standard stochastic optimization frameworks, including MDPs and stochastic programming. Nevertheless, estimated RLDs can be biased due to multiple factors, such as noisy sensors or sparsity of training data. This can result in that the decisions computed using  MDPs or stochastic programming can be far from true optimal. 2) The existing works deal with multiple identical components, thus needing just one spare inventory. However, most industrial systems have multiple \textit{non-identical} components, so they need to manage multiple spare inventories accounting for the diversity of spare parts and their potential dependencies. This paper addresses these two critical open issues.

\subsection{Overview and Contributions}
This paper develops a Distributionally Robust Chance Constrained (DRCC) formulation to jointly optimize CBM schedules and spare parts provisioning. The uncertainty of this formulation is associated with the remaining lifetime of components.  This uncertainty is characterized by prognostic models with an online updating framework that utilizes real-time CM data to periodically estimate components’ RLDs. The optimization model can be implemented in a rolling horizon fashion to adapt its decisions to the degradation process experienced by the components.

In contrast to previous studies, our proposed DRCC formulation acknowledges that estimated RLDs can be \textit{biased} and so seeks optimal solutions that are robust against distributional perturbations within a data-driven ambiguity set. Moreover, the formulation deals with a system with multiple machines composed of  \textit{non-identical} components. Hence, the resulting model allows the coordination of multiple spare part inventories, accounting for the diversity of components in real-world industrial applications. We capture the maintenance cost dependence between the different constituting components of a machine by adopting an opportunistic approach, which schedules CBM by incentivizing the grouping of repair activities when profitable, i.e., repairing more than one component when a machine is shut down.

Our problem formulation includes two Distributionally Robust (DR) chance constraints that aim to restrict the probability of the following undesired events: (1) \edit{downtime} of critical components and (2) unexpected failures. We prove that the DR chance constraint associated with the \edit{downtime} of critical components admits an exact mixed-integer linear programming representation. We demonstrate that the DR chance constraint related to unexpected failures can be rewritten as a regular chance constraint whose underlying distribution is computed by solving a \edit{Distributionally Robust Optimization (DRO)} problem. We provide an explicit formula to calculate that distribution. Combining these results, we show that the proposed DRCC problem admits a MILP reformulation that can be efficiently solved by off-the-shelf optimizer software. Notice that this MILP reformulation is derived for general discrepancy-based ambiguity sets. Finally, for computational experiments, we adopt a type-$\infty$ Wasserstein ambiguity set and derive closed-form expressions to the parameters of the MILP reformulation. The efficacy of our methodology is showcased in a \edit{simulation} case study on wind turbines, where we analyze the impact of using the proposed DRCC formulation versus other benchmark models based on stochastic programming and robust optimization. 

In summary, the main contributions of this paper are threefold:
\begin{itemize} 
\item This paper develops a DRCC formulation to jointly optimize CBM schedules and spare parts provisioning, integrating insights from predictive analytics. The proposed DRCC formulation recognizes the potential \textit{bias} in RLDs estimated by predictive analytics. Consequently, it seeks maintenance and spare parts decisions robust against distribution perturbations within discrepancy-based ambiguity sets. This acknowledgment of bias contributes to the formulation's robustness, addressing a critical aspect overlooked in previous studies. \edit{Additionally, the proposed formulation can be applied to purely CBM scheduling, offering a methodology for determining DR maintenance schedules that incorporate predictive analytics insights. }
\item The formulation tackles the complexity of systems with multiple machines, each comprised of \textit{non-identical} components. This feature distinguishes our model by facilitating the coordination of diverse spare part inventories.  Additionally, it captures the interdependence of maintenance costs among the various components within a machine. This is achieved through an opportunistic approach that strategically schedules CBM, encouraging the simultaneous repair of multiple components during machine shutdowns when economically convenient.
\item The paper shows that the proposed formulation admits an exact MILP reformulation. We provide explicit formulas to implement the optimization model for general discrepancy-based ambiguity sets. Interestingly, we demonstrate that the DR related to unexpected failures can be rewritten as a regular chance constraint whose underlying distribution can be computed explicitly by solving a  simple DRO problem. Finally, for the numerical illustration, we study a type-$\infty$ Wasserstein ambiguity set and derive closed-form expressions for the parameters of the MILP reformulation. \edit{The model is applied to a case study on wind turbines constructed with simulated degradation data.} 
\end{itemize}

The rest of the paper is organized as follows. Section \ref{sec:prob_statement} details the \edit{problem setting} for the joint optimization model. Section \ref{sec:formulation} introduces the proposed DRCC formulation for jointly optimizing CBM and spare provisioning. Section \ref{sec:computational_studies} discusses the results obtained in multiple simulation studies. Finally, Section \ref{sec:conclusions} concludes this paper with some closing remarks. 
\section{Problem Setting}\label{sec:prob_statement}

\begin{figure}[htbp]
    \centering
    \includegraphics[width=0.7\linewidth]{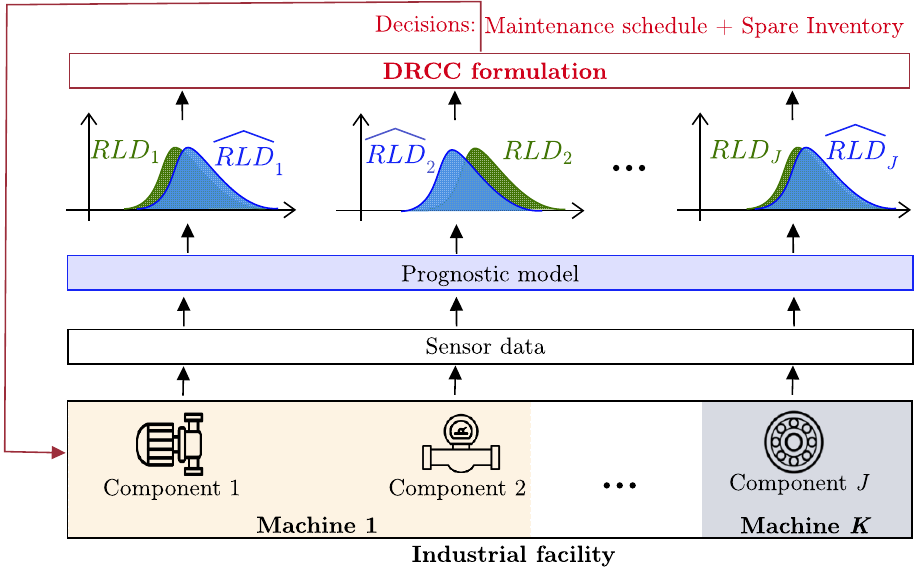}
    \caption{Graphical summary of the problem setting.}
    \label{fig:proposal}
\end{figure}

We consider a problem setting where the objective is to jointly optimize CBM and spare parts provisioning for a set of $K$ machines,  indexed by $k \in \mathcal{K}=\{1,...,K\}$,  \edit{over a finite planning horizon $\mathcal{T}=\{1,...,T_{max}\}$}. Each machine consists of different critical components (e.g., bearings, valves, pumps, etc.) that are monitored by sensors. The subset of components that constitutes the machine $k$ is denoted as $\mathcal{J}_k,\; k \in \mathcal{K}.$ \edit{It is assumed that the machines do not share installed components (i.e., an installed component cannot simultaneously contribute to the functionality of two different machines); however, they may share spare parts (i.e., the same spare part can be used to repair two different machines).} Thus, the collection of subsets $\mathcal{J}_1,...,\mathcal{J}_K$ satisfies $\mathcal{J}_{k_1} \cap \mathcal{J}_{k_2}  = \emptyset,\; k_1 \neq k_2 \in \mathcal{K}.$ The set of all components under study is indexed by $j \in \mathcal{J}=\{1,...,J\}$  and can be determined by $\mathcal{J}=\bigcup_{k \in \mathcal{K}}\mathcal{J}_k$. The problem setting is summarized in Figure \ref{fig:proposal}, and the building blocks and assumptions are discussed next.

 We assume that components undergo gradual degradation over time. Condition monitoring data from sensors attached to the components exhibit trends that are correlated with the severity of the underlying physical degradation process of the component that is being monitored. Given a predefined failure threshold, trends in the CM data can be used to predict the RLD by prognostic models as described later in Section \ref{subsec:prognostics}. The remaining lifetime of component $j \in \mathcal{J}$ is modeled as a non-negative random variable, denoted as $\tilde{\omega}_j$ with support $[0,\tau_j^{\text{max}}]=\Xi_j \subseteq \mathbb{R}_+$. We let $\tilde{\Omega} =[\tilde{\omega}_1,...,\tilde{\omega}_j]$ be a random vector with support $\Xi = \Xi_1 \times \Xi_2 \times \cdots \times \Xi_J$ representing the remaining lifetimes for all components. We assume that the random variables $ \tilde{\omega}_1,...,\tilde{\omega}_j$ are independent, which is a common assumption in maintenance optimization problems (e.g., \cite{keizer2017joint,zhang2017joint}).

 For a given component $j \in \mathcal{J}$, the corresponding prognostic model computes $\mathrm{\widehat{RLD}_j}$ that is a data-driven estimate of the ground truth distribution, $\mathrm{RLD}_j$, which governs the uncertainty of the component remaining lifetime $\tilde{\omega}_j$.

 Estimating the uncertainty of components' remaining lifetime is key to optimizing maintenance schedules. It enables modeling the trade-off between conservative and risk-taking maintenance policies. Conservative policies tend to repair components earlier than needed, thereby reducing the risk of unexpected failures but losing partial useful life of components. Risk-taking policies, on the other hand, postpone components' repairs as much as possible. This results in a more efficient usage of the entire lifetime of components but in a higher risk of unexpected failure. Maintenance of individual components involves shutting down the entire machine, which can sometimes be costly. Consequently, we adopt an opportunistic maintenance policy that encourages the replacement of partially degraded components when there is a convenient opportunity to do so. For instance, when a machine is idle due to the repair of a degraded component, technicians can leverage the opportunity to replace other partially degraded components within the  machine. These opportunistic maintenance decisions are determined by the optimization model.

For the spare part inventory, we assume that some machines may need the same spare parts (\edit{i.e., they may share spare parts}), so we do not need to manage $J$ different spare inventories. Instead, we analyze $L$ different spare inventories accounting for the $L$ different types of components distributed across the different machines. The subset of components of type $l$ is denoted as $\mathcal{M}_l \subseteq \mathcal{J}, \; l \in \mathcal{L}=[1,...,L]$. There are holding costs associated with storing spares because inventory managers are interested in just-in-time orders, i.e., in aligning the arrival of spare orders with the maintenance schedule. The spare inventories can be replenished by two different types of spare orders, referred to as regular orders and expedited orders. Regular orders are cheaper per unit than expedited orders, yet have larger lead times.

 Our final goal is to integrate the RLDs estimated by prognostic models (i.e., $\mathrm{\widehat{RLD}}_j, j\in \mathcal{J}$ ) into a decision framework built to jointly optimize maintenance and spare inventory. By doing this, decisions on the maintenance schedule (when to repair) and spare inventory (when to place regular and expedited orders) can be adapted according to the insights obtained from prognostic models. In contrast to previous studies, our problem setup acknowledges that the distribution estimated by the prognostic algorithm can be biased. Hence, optimization techniques such as stochastic programming or simulation-based methods \textit{cannot} guarantee good performance in this setup. To tackle this issue, we propose using a DRO formulation with a finite planning horizon. \edit{ This formulation is introduced in Section \ref{sec:formulation} and relies on two important assumptions, which are listed and discussed next.
   \begin{assum}\label{as:one_failure}
If a component has been repaired within the planning horizon, it cannot fail within the remaining periods of the same planning horizon.
\end{assum}
\begin{assum}\label{as:one_repair}
  Each component needs a single repair within the planning horizon.
 \end{assum}
These two assumptions are common in the existing literature on maintenance optimization problems (e.g., \cite{shi2022stochastic, bakir2021integrated}). Although they may seem restrictive, they can be justified in practical applications by combining a careful selection of the planning horizon length with the execution of the decision model in a rolling horizon fashion. More specifically, it is sufficient to choose a planning horizon shorter than the minimum guaranteed lifetime of the components under study. This ensures that if a component is replaced within a given planning horizon, it cannot fail again within the remaining periods of the same planning horizon (Assumption \ref{as:one_failure}) and that each component will require at most one repair within the planning horizon (similar to Assumption \ref{as:one_repair}). However, under Assumption \ref{as:one_repair}, components in good condition are enforced to receive one repair within the planning horizon. This issue is addressed by implementing the optimization strategy in a rolling-horizon framework with a fixed freeze period. In this approach, only decisions committed within the freeze period are executed, while decisions beyond the freeze period are updated as the freeze period progresses. For components in good condition, maintenance activities will be scheduled beyond the freeze period. Therefore, these maintenance activities will never be performed in practice.
}


 \subsection{Sensor-driven Failure Prognostic Models}\label{subsec:prognostics}

 Many industrial components undergo gradual degradation over time. It is often possible to track the evolution of the degradation process by modeling the associated CM data. CM data usually exhibit patterns and trends that are correlated with the severity of the underlying physical degradation process. \add{Prognostic modeling approaches leverage these trends}\remove{A prognostic model can be utilized} to estimate the RLDs of partially degraded components. \add{The joint optimization model developed in this paper} admits any prognostic modeling framework capable of predicting RLDs. 
 
 There are many classes of prognostic models that can be used to predict RLDs. They rely on various methodologies ranging from stochastic methods and statistical approaches to Machine Learning tools and physics-based models--A survey on this topic is provided in \citep{elattar2016prognostics}. For expository purposes, we focus on the stochastic prognostic modeling frameworks. These models assume that CM data evolves according to some stochastic process that reflects the uncertainty of the underlying physical degradation process. In other words, the degradation process of a component $j \in \mathcal{J}$ is modeled by a random process $\zeta(t,\Theta_j)$, where $t$ is the cumulative operating time of the component, and $\Theta_j$ is a set of stochastic model parameters that follow a predefined prior distribution estimated from historical CM data. Most of the stochastic prognostic models have an integrated online updating mechanism that is usually based on a Bayesian framework. In this updating framework, the prior distribution of $\Theta_j$ is periodically updated based on recently observed CM data, $\mathcal{O}_j^{t}$, where $\mathcal{O}_j^{t}$ is the set of CM data collected from component $j$ up to time $t$.  The online updating process improves the accuracy of the predicted RLD because it incorporates data associated with the most recent degradation state of the component that is being monitored. This gives the flexibility to revise maintenance and spare parts ordering decisions based on the latest observed degradation states. In other words, each time the RLDs are updated, the joint optimization model is resolved to generate revised decisions.

  Failure is assumed to occur when $\zeta(t,\Theta_j)$ crosses a predefined failure threshold, $\Lambda_j$. In most applications, this threshold is defined in terms of the maximum allowable amplitude of the CM data. The RLD of component $j$ can be calculated as the first-passage-time event. Specifically, the remaining lifetime, $\tilde{\omega}_j$, can be expressed as $\tilde{\omega}_j = \min_{t \geq 0} {\zeta(t, \Theta_j) \geq \Lambda_j}$. The RLD updated at time $t$ is evaluated by conditioning on the CM data $\mathcal{O}_j^{t}$ as follows:
\begin{equation}\label{eq:prognostic}
\begin{split}
\widehat{RLD}_j = \mathbb{P}(\tilde{\omega}_j \leq \tau| \mathcal{O}_j^{t})&= 1-\mathbb{P}(\tilde{\omega}_j > \tau| \mathcal{O}_j^{t})\\
&= 1 - \int \mathbb{P}( \zeta(s,\Theta_j) \leq \Lambda_j,\; 0\leq s \leq \tau )  \mathbb{P}(\Theta_j|\mathcal{O}_j^{t}) d \Theta_j.
\end{split}
\end{equation}
We highlight the flexibility of our optimization framework by noting that multiple prognostic modeling frameworks can be used to solve \eqref{eq:prognostic} and compute the RLD. In particular, well-established prognostic methodologies based on Brownian motion \citep{gebraeel2005residual} and Particle Filter \citep{orchard2009particle} can be directly integrated into our optimization framework. This allows us to model the degradation of different types of components, thereby enabling the optimization of maintenance and spare provisioning for a large variety of industrial components.


\section{A Distributionally Robust Chance Constrained  Formulation}\label{sec:formulation}

In this section, we present our Distributionally Robust Chance Constrained (DRCC)  formulation to jointly optimize CBM schedules and spare parts provisioning. We begin by defining the ambiguity set used to define model robustness. We then discuss the objective function of our optimization model followed by the model constraints.

\subsection{Ambiguity Set}\label{subsec:ambiguityset}

The uncertainty of our optimization model is associated with the remaining lifetime of components, which is modeled as a vector of non-negative independent random variables
$\tilde{\Omega} =(\tilde{\omega}_1,...,\tilde{\omega}_J)^\top$. For any component $j \in \mathcal{J}$, the distribution of $\tilde{\omega}_j$ is unknown but can be estimated from sensor data by a prognostic model, resulting in $\widehat{\mathrm{RLD}}_j$ as described in Section \ref{subsec:prognostics}. The joint distribution of $\tilde{\Omega}$, denoted as $\widehat{\mathrm{RLD}}$, can be then estimated by using the independence of $\tilde{\omega}_j, \; j \in \mathcal{J}$.

To model potential estimation errors in the distribution $\widehat{\mathrm{RLD}}$, we construct an ambiguity set to capture possible perturbations of this distribution. To this end, we rely on the independence of the RLD of each component to propose an ambiguity set with the following structure:
\begin{equation}\label{eq:global_ambiguity_set}
    \mathcal{P} = \left \{ \mathbb{P}:\; \mathbb{P}\left \{\tilde{\Omega} \in \Xi_1\otimes \Xi_2 \otimes  \cdots \otimes \Xi_J\right\}=1, \;  \mathbb{P} \in  \mathcal{P}_1 \otimes \mathcal{P}_2 \otimes \cdots \otimes \mathcal{P}_J  \right\},
\end{equation}
where
\begin{equation}\label{eq:local_ambiguity_set}
   \mathcal{P}_j = \left \{ \mathbb{P}: \; \mathbb{P}\left \{\tilde{\omega}_j \in \Xi_j \right\}=1,\; \varphi(\mathbb{P},\widehat{\mathrm{RLD}}_j) \leq \delta_j \right\}, \; j \in \mathcal{J}.  
\end{equation}
Notice that $\mathcal{P}_j$ corresponds to the ambiguity set associated with the random variable $\tilde{\omega}_j$. $\mathcal{P}_j$ is defined using the functional $\varphi(,)$ that measures the discrepancy between any candidate distribution $\mathbb{P}$ and $\widehat{\mathrm{RLD}}_j$. Thus, $\mathcal{P}_j$ includes any distribution $\mathbb{P}$ with support $\Xi_j$ such that its discrepancy with respect to the estimated distribution $\widehat{\mathrm{RLD}}_j$ is at most  $\delta_j \geq 0$. The functional $\varphi(,)$ can be defined using different approaches, including total variation metric, Wasserstein distance, $\phi-$divergence distance, or distribution moments \citep{rahimian2019distributionally}. The hyper-parameter $\delta_j $ controls the size of the ambiguity set.  It can be used to represent our confidence in the prognostic model that is being used to predict the RLD of component $j$. A small $\delta_j$ will correspond to a high level of confidence in our prognostic results (and vice versa).   

\subsection{Decision Variables and Objective Function}

In this subsection, we detail both the decision variables and the objective function considered in our formulation. Let $\mathcal{T}$ be the set of maintenance epochs. First, we specify decision variables concerning maintenance schedule. We introduce the binary decision variable $x_{jt}$ that is one if component $j \in \mathcal{J}$ is scheduled to be repaired at time $t \in \mathcal{T}$. To repair component $j$, it is necessary to shut down the machine that contains component $j$. Thus, we define the binary decision variable $y_{kt}$ that is one if machine $k \in \mathcal{K}$ is shut down at time $t \in \mathcal{T}$. This decision variable is utilized to model the machine-component dependencies and to encourage opportunistic repairs when profitable, i.e., repairing more than one component when a machine is shut down. We define a binary decision variable $z_{t}$ that is one if the maintenance crew is deployed for conducting repairs at decision epoch $t \in \mathcal{T}$. 
This decision variable is used to capture the deployment or setup cost of the maintenance crew, which is relevant to encourage group repairs of different machines.  Secondly, we introduce decision variables related to spare parts management. We define the integer decision variable $h_{lt}, \; l \in \mathcal{L}, \; t \in \mathcal{T}$ that captures the inventory level of component type $l$ at decision epoch $t$. We introduce the integer decision variable $g^{reg}_{lt} (g^{exp}_{lt}), \; l \in \mathcal{L}, \; t \in \mathcal{T}$ that represents the regular (expedited) spare part order of component type $l$ placed at decision epoch $t$. Finally, we define the decision variable $r_t$ that is one if a regular spare order is placed at $t \in \mathcal{T}$. This variable will be used to capture the fixed costs associated with regular orders. The main notation used in this formulation is summarized by Table \ref{tab:notation}.
 {
\renewcommand{\arraystretch}{0.95}
\begin{table}[htbp]
\caption{Notation}
\label{tab:notation}
\centering 
\begin{tabular}{ll}
\hline
Notation & Description \\ \hline
\text{\textit{Sets}:}& \vspace{0.1cm}\\
$\mathcal{J}$ & : Set of components.\\
$\mathcal{K}$ & : Set of machines.\\
$\mathcal{L}$ & : Set of type of spare parts.\\
$\mathcal{T}$ & : Set of maintenance epochs. \\
$\mathcal{N}$ & : Set of samples used by $\widehat{\mathbb{P}}_N (\tilde{\Omega})$.\vspace{0.1cm}\\
\text{\textit{Decision variables}:}& \vspace{0.1cm} \\
$x_{jt} \in \{0,1\}$ & : 1, if component $j$ undergoes repair at time $t$. \\
$y_{kt} \in \{0,1\}$ & : 1, if machine $k$ is shut down at time $t$. \\
$z_{t} \in \{0,1\}$ & : 1, if crew is deployed for maintenance at time $t$. \\ 
$h_{lt} \in \mathbb{Z}_{+}$ & : Inventory level of spare type $l$ at time $t$. \\
$g^{reg}_{lt} \in \mathbb{Z}_{+}$ & : Regular order of spare type $l$ placed at time $t$. \\
$g^{exp}_{lt} \in \mathbb{Z}_{+}$ & : Expedited order of spare type $l$ at time $t$. \\
$r_t \in  \{0,1\}$ &: 1, if some regular spare order is placed at time $t$. \\
$a_{je} \in  [0,1]$ &: Auxiliary variable used to model the second chance constraint. \\ 
$\bm{X}$ & : Vector of all decision variables.\\
\text{\textit{Parameters}:}& \vspace{0.1cm} \\
$J$ &: Total number of components under study. \\
$K$ &: Total number of machines under study. \\
$N$ &: Total number of samples. \\
$T_{max}$ &: Length of the planning horizon. \\
$\delta_j$ &: Parameter controlling the size of the ambiguity set. \\
$\Delta^{upd}$ &: Length of the freeze period. \\
$\Delta^{reg}$ &: Lead time for regular orders. \\
$C^{pr}_j$ &: Preventive maintenance cost of component $j$. \\
$V^{pr}_j$ &: Cost of premature repair for component $j$. \\
$C^{co}_j$ &: Corrective maintenance cost of component $j$. \\
$V^{co}_j$ &: Cost of late repair for component $j$. \\
$C^{down}_k$ &: Cost of shutting down machine $k$ for maintenance. \\
$C^{crew}$ &: Setup cost. \\
$C^{hold}_l$ &: Holding cost for component type $l$. \\
$C^{reg}_l$ &: Per-unit cost of regular orders for component type $l$. \\
$C^{exp}_l$ &: Per-unit cost of expedited orders for component type $l$. \\
$B^{reg}$ &: Fixed cost for placing regular spare part orders. \\
$M$ &: Maximum capacity of the maintenance crew. \\
$G_t$ &: Maximum capacity of spare part suppliers at decision epoch $t$. \\
$\rho$ &: Threshold on the maximum \edit{downtime} of each component. \\
$\gamma$ &: Threshold on the maximum number of unexpected failures. \\
$\epsilon$ &: Confidence level for the first chance constraint. \\
$\beta$ &: Confidence level for the second chance constraint.  
\vspace{0.1cm}  \\
   \hline
\end{tabular}
\end{table} 
}

We aim to minimize the operational costs associated with maintenance tasks and spare parts management. For this, we propose the cost function presented in \eqref{eq:obj1}, where $\bm{X}$ is a vector that contains all decision variables \edit{and $\mathbb{I}(\cdot)$ represents the indicator function}. The cost function is divided into six parts, which are described next.

\noindent(i) \textit{Preventive maintenance cost} captures the cost of preventive repairs, which are repairs carried out before failure, i.e., $t < \tilde{\omega}_j$. This cost consists of a fixed penalty $C_j^{pr}$ plus a variable penalty $V_j^{pr}  (\tilde{\omega}_j-t)$ that depends on how early the repair takes place. \edit{This variable penalty serves as a mechanism to discourage overly conservative maintenance decisions.}

\noindent(ii) \textit{Corrective maintenance cost} is incurred when repairs are performed in response to unexpected failures, i.e., $t \geq \tilde{\omega}_j$. This cost also consists of a fixed cost $C_j^{co}$ plus a variable penalty $V_j^{co}  (t- \tilde{\omega}_j)$ that depends on how late the repair takes place. In practice, we usually see $C_j^{pr}\leq C_j^{co}$ and $V_j^{pr} \leq V_j^{co}$ because unexpected failures often take longer to repair and may cause secondary failures. Notice that preventive and corrective maintenance costs model the trade-off between early and late repairs.

\noindent(iii) \textit{Machine shutdown cost} ($C^{down}_k$) is incurred when a machine is shut down to conduct maintenance to any of its components. This cost captures the machine-component dependence and gives an incentive to conduct opportunistic repairs to some components.

\noindent(iv) \textit{Setup Cost} models the fixed cost ($C^{crew}$) incurred for conducting maintenance. This cost is used to encourage group repairs at the machine level.

\noindent(v) \textit{Holding cost} accounts for the cost ($C_l^{hol}$) incurred for keeping spare parts inventory.

\noindent(vi) \textit{Spare part order cost} captures the cost of regular ($C_l^{reg}$) and expedited orders ($C_l^{exp}$). \edit{In practice, we typically observe that 
$C_l^{\text{reg}} \leq C_l^{\text{exp}},$ 
as expedited orders require rapid delivery, significantly increasing shipping costs. Consequently, expedited orders are placed only when an unpostponable repair is needed, and the required spare part is unavailable in on-hand inventory. If the repair is not pressing to satisfy reliability requirements, the optimization model will postpone the repair, waiting for a regular spare part order to arrive and perform the maintenance activity. Regular orders are less expensive than expedited ones, but suppliers often impose a fixed cost, 
$B^{\text{reg}},$ 
on regular orders to encourage grouping---i.e., ordering multiple units in a single transaction. This creates a trade-off between holding inventory and placing orders more frequently. Expedited orders, on the other hand, do not incur this fixed cost because the per-unit cost is high, encouraging the inventory manager to purchase no more than one unit per expedited order.

}

\begin{equation}\label{eq:obj1}
\begin{split}
   \mkern-10mu f(\bm{X},\tilde{\Omega})  = & \overbrace{\sum_{j \in \mathcal{J}}\sum_{t \in \mathcal{T}} \left( \left( V^{pr}_j  (\tilde{\omega}_j-t) + C^{pr}_j \right)  \mathbb{I}\{ \tilde{\omega}_j >  t\} \right)   x_{jt} }^{\text{(i) Preventive maintenance cost}}+ \\ &  \overbrace{\sum_{j \in \mathcal{J}}\sum_{t \in \mathcal{T}} \left(  \left(  V^{co}_j  (t-\tilde{\omega}_j) +C^{co}_j  \right) \mathbb{I}\{ \tilde{\omega}_j \leq  t\}  \right)   x_{jt} }^{\text{(ii) Corrective maintenance cost}} +\mkern-18mu \overbrace{\sum_{k \in \mathcal{K}}\sum_{t \in \mathcal{T}} C^{down}_k  y_{kt}}^{\text{(iii) Machine shutdown cost}} +\overbrace{\sum_{t \in \mathcal{T}} C^{crew}  z_t}^{\text{(iv) Setup cost}}+ \\
    &\overbrace{\sum_{l\in \mathcal{L}}\sum_{t \in \mathcal{T}} C_l^{hol} h_{lt}}^{\text{(v) Holding cost}} 
    +\overbrace{ \sum_{l\in \mathcal{L}}\sum_{t \in \mathcal{T}} \left (C_l^{reg} g^{reg}_{lt} + C_l^{exp} g^{exp}_{lt} \right )  +\sum_{t \in \mathcal{T}}B^{reg} r_t }^{\text{(vi) Spare part order cost}}.
\end{split}
\end{equation}
Note that $f(\bm{X},\tilde{\Omega})$ can be rewritten as follows:
\begin{equation*}\label{eq:obj_2}
\begin{split}
    f(\bm{X},\tilde{\Omega})  = & \sum_{j \in \mathcal{J}}\sum_{t \in \mathcal{T}}  \alpha(\tilde{\omega}_j,t)  x_{jt}   +   \sum_{k \in \mathcal{K}}\sum_{t \in \mathcal{T}} C^{down}_k  y_{kt}  + \sum_{t \in \mathcal{T}} \left ( C^{crew}  z_t +B^{reg} r_t  \right )+ \\ & \sum_{l\in \mathcal{L}}\sum_{t \in \mathcal{T}} \left (C_l^{hol} h_{lt} + C_l^{reg} g^{reg}_{lt} + C_l^{exp} g^{exp}_{lt} \right ),   
\end{split}
\end{equation*}
where
\begin{equation*} 
\begin{split}
    \alpha(\tilde{\omega}_j,t)=  \left( V^{pr}_j(\tilde{\omega}_j-t) + C^{pr}_j \right) \mathbb{I}\{ \tilde{\omega}_j >  t\} + \left(  V^{co}_j(t-\tilde{\omega}_j) +C^{co}_j  \right) \mathbb{I}\{ \tilde{\omega}_j \leq  t\}.
\end{split}
\end{equation*}
The cost function $ f(\bm{X},\tilde{\Omega})$ is stochastic due to its dependence on $\tilde{\Omega}$, which influences the preventive and corrective maintenance costs. Then, we aim to solve the following DR optimization problem: 

\begin{equation}\label{eq:objective}
\min_{\bm{X} \in \mathcal{X}} \left\{ \sup_{\mathbb{P} \in \mathcal{P} }{\mathbb{E}_{\mathbb{P}}[f(\bm{X},\tilde{\Omega})]} \right\},
\end{equation}
where $\mathcal{X}\subseteq \mathbb{R}^{|\bm{X}|}$ corresponds to the feasible region of $\bm{X}$ determined by the constraints explained later. We now focus on deriving a tractable expression for $\sup_{\mathbb{P} \in \mathcal{P} }{\mathbb{E}_{\mathbb{P}}[f(\bm{X},\tilde{\Omega})]}$. As the uncertainty of $f(\bm{X},\tilde{\Omega})$ is solely associated with $\alpha(\tilde{\omega}_j,t)$, the following equality holds.
   \begin{equation}\label{eq:sup1}
   \begin{split}
     \sup_{\mathbb{P} \in \mathcal{P} }  \mathbb{E}_{\mathbb{P}}[f(\bm{X},\tilde{\Omega})]=& \sum_{k \in \mathcal{K}}\sum_{t \in \mathcal{T}} C^{down}_k  y_{kt}  + \sum_{t \in \mathcal{T}} \left ( C^{crew}  z_t +B^{reg} r_t  \right )+ \\ & \sum_{l\in \mathcal{L}}\sum_{t \in \mathcal{T}} \left (C_l^{hol} h_{lt} + C_l^{reg} g^{reg}_{lt} + C_l^{exp} g^{exp}_{lt} \right )    + \\ & \sup_{\mathbb{P} \in \mathcal{P} }{\mathbb{E}_{\mathbb{P}}  \left[\sum_{j \in \mathcal{J}}\sum_{t \in \mathcal{T}} \alpha(\tilde{\omega}_j,t)  x_{jt}\right]. }
     \end{split}
\end{equation} 
\edit{Most of the terms on the right-hand side of \eqref{eq:sup1} are linear, except for the term $\sup (\cdot)$. However, this supremum can be simplified to a linear combination of the decision variables with pre-computed weights, as shown next in Proposition \ref{prop:obj}. }
\begin{prop}\label{prop:obj}
 If $\sum_{t \in \mathcal{T}}  x_{jt}=1$ for any $j \in \mathcal{J}$, then
\begin{equation*}
      \sup_{\mathbb{P} \in \mathcal{P} }{\mathbb{E}_{\mathbb{P}}  \left[\sum_{j \in \mathcal{J}}\sum_{t \in \mathcal{T}} \alpha(\tilde{\omega}_j,t)  x_{jt}\right] }=   \sum_{j \in \mathcal{J}}  \sum_{t \in \mathcal{T}}  \psi_{jt} x_{jt},
\end{equation*} 
where $ \psi_{jt}=\sup_{\mathbb{P} \in \mathcal{P}_j}{\mathbb{E}_{\mathbb{P}}   \left[\alpha(\tilde{\omega}_j,t) \right]}, \; j \in \mathcal{J}, \; t \in \mathcal{T}$  are constant values that can be pre-computed.
\end{prop}
 
\begin{proof}
Note that
\begin{subequations}
\begin{align} 
\sup_{\mathbb{P} \in \mathcal{P}}{\mathbb{E}_{\mathbb{P}}  \left[\sum_{j \in \mathcal{J}}\sum_{t \in \mathcal{T}} \alpha(\tilde{\omega}_j,t)  x_{jt}\right] }   & = \sum_{j \in \mathcal{J}}  \sup_{\mathbb{P} \in \mathcal{P}_j}{\mathbb{E}_{\mathbb{P}}   \left[\sum_{t \in \mathcal{T}} \alpha(\tilde{\omega}_j,t)  x_{jt} \right] } \label{eq:o1}\\
& = \sum_{j \in \mathcal{J}} \sum_{t \in \mathcal{T}}  \sup_{\mathbb{P} \in \mathcal{P}_j}{\mathbb{E}_{\mathbb{P}}   \left[\alpha(\tilde{\omega}_j,t) \right]  x_{jt} },  \label{eq:o2} 
\end{align}
\end{subequations}
where \eqref{eq:o1} follows from the definition of using  $\mathcal{P}$ in \eqref{eq:global_ambiguity_set}. Then, as $\sum_{t \in \mathcal{T}}  x_{jt}=1$, only one $x_{jt}$ is going to be equal 1, so we can then determine the supremum for each time $t$, resulting in \eqref{eq:o2}. The proof finishes by defining $\psi_{jt}=\sup_{\mathbb{P} \in \mathcal{P}_j}{\mathbb{E}_{\mathbb{P}}   \left[\alpha(\tilde{\omega}_j,t) \right]}.$ 
\end{proof}
\edit{The weights $\psi_{jt}, \; j \in \mathcal{J}, \; t \in \mathcal{T}$, represent the worst-case the expected maintenance costs and can be determined by solving a simple DRO problem. The condition $\sum_{t \in \mathcal{T}} x_{jt} = 1$ for any $j \in \mathcal{J}$ stated in Proposition \ref{prop:obj} implies that each component must be repaired exactly once within the planning horizon. This condition is consistent with Assumption \ref{as:one_repair}}. Therefore, we can invoke Proposition \ref{prop:obj} to easily show that the optimization problem of \eqref{eq:objective} reduces to:
{
\begin{multline*} 
 \min_{\bm{X} \in \mathcal{X}} \left\{ \sum_{k \in \mathcal{K}}\sum_{t \in \mathcal{T}} C^{down}_k  y_{kt}  + \sum_{t \in \mathcal{T}} \left ( C^{crew}  z_t +B^{reg} r_t  \right )+ \sum_{l\in \mathcal{L}}\sum_{t \in \mathcal{T}} \left (C_l^{hol} h_{lt} + C_l^{reg} g^{reg}_{lt} + C_l^{exp} g^{exp}_{lt} \right )  \right.  +\\   \left. \sum_{j \in \mathcal{J}}\sum_{t \in \mathcal{T}}   \psi_{jt}  x_{jt}  \right\},
\end{multline*}
which corresponds to a linear objective function minimized over $\mathcal{X}$. Next, we describe the constraints used to describe feasible region $\mathcal{X}$.

\subsection{Deterministic Constraints }

We introduce the deterministic constraints considered in this formulation. For notational convenience, the feasible region induced by these deterministic constraints is denoted as $\mathcal{D}$.

We first enforce that any component needs exactly one repair within the planning horizon. This is modeled by \eqref{eq:onerepair}. 
\begin{equation}\label{eq:onerepair}
    \sum_{t \in \mathcal{T}} x_{jt}=1, \quad j \in \mathcal{J}.
\end{equation}
Repair tasks are scheduled at time $t$ only if the maintenance crew is deployed but the maintenance crew can conduct at most $M$ repairs per maintenance epoch. These two facts are captured by \eqref{eq:capacity}. 
\begin{equation}\label{eq:capacity}
    \sum_{j \in \mathcal{J}} x_{jt} \leq M  z_t  , \quad t \in \mathcal{T}.
\end{equation}
We also model the component-machine dependence by including constraint \eqref{eq:component-machine}. This captures the fact that the machine must be shut down in order to schedule any repair of a machine component.

\begin{equation}\label{eq:component-machine}
    \sum_{j \in \mathcal{J}_k}x_{jt} \leq |\mathcal{J}_k|  y_{kt}, \quad  k \in \mathcal{K}, \; t \in \mathcal{T}.
\end{equation}

\edit{Additional constraints are required for spare parts inventory. Specifically, the total number of spare parts of a regular order placed at time $t$ cannot exceed the supplier's maximum capacity, $G_t, \; t \in \mathcal{T}$. Moreover, a regular order can only be placed if $r_t = 1$, where $r_t$ is the binary decision variable representing the payment of the fixed order cost for regular orders at time $t$. These two conditions are enforced by \eqref{eq:fixed_order}.  We characterize the evolution of the inventory level over time by constraint \eqref{eq:inventory}, which includes the possibility of regular and expedited spare part orders.  Our formulation assumes that the lead times are deterministic, fixed, and known. Regular spare orders have a lead time denoted by $\Delta^{reg}>0$, while expedited orders are assumed to be instantaneous, with a lead time of zero. Orders always arrive at the beginning of each period, so, if needed, they are available for installation within the same period. Mention that $\Tilde{g}_{l,t-\Delta^{reg}}^{reg}, \;  1\leq t \leq \Delta^{reg} $ correspond to regular spare parts orders placed in the past, which cannot be modified and will be arriving at time $t$.   It is important to note that real-world applications often involve stochastic lead times, which pose additional challenges in coordinating maintenance and spare inventory. This important issue is out of the scope of this paper but represents an interesting direction for future extensions of this work.

\begin{equation}\label{eq:fixed_order}
    \sum_{l \in \mathcal{L}} g_{lt}^{reg} \leq G_t r_t, \quad t \in \mathcal{T}.
\end{equation}

\begin{equation}\label{eq:inventory}
    h_{lt}=\begin{cases} 
      h_{l0} +g_{l,1}^{exp}+\Tilde{g}_{l,t-\Delta^{reg}}^{reg}- \sum_{j \in \mathcal{M}_l}x_{j,1},& 1\leq t \leq \Delta^{reg}   \\
      h_{l,t-1} +g_{lt}^{exp}+g_{l,t-\Delta^{reg}}^{reg}- \sum_{j \in \mathcal{M}_l}x_{jt}, &  \Delta^{reg} +1\leq t   
   \end{cases},\; l\in \mathcal{L}, \; t\in \mathcal{T}.
\end{equation}}

 \subsection{Distributionally Robust Chance Constraints}
 
Now we detail the distributionally robust chance constraints considered in our formulation.

 \subsubsection{Limiting the \edit{Downtime} of Individual Components:}

When a component suffers an unexpected failure, the underlying machine must remain shut down until the component is repaired. Therefore, the \edit{downtime} of one component impacts the availability of the machine, negatively affecting the production processes. To restrict the probability of long \edit{downtime} of components due to \edit{an unexpected failure}, we introduce the following DRCC set:
\begin{equation*}
\begin{split}
    \mathcal{Z}_1= \Biggl\{\bm{X} \in \mathcal{D}: \inf_{\mathbb{P} \in \mathcal{P} } \mathbb{P}\biggl\{ \tilde{\Omega}: \;  \sum_{t \in \mathcal{T}}&[t-\tilde{\omega}_j]_{+}  x_{jt}\leq  \rho \biggr\}  \geq 1-\epsilon, \quad j \in \mathcal{J} \Biggr\},
    \end{split}
\end{equation*}
\edit{where $[\cdot]_+= \max (0,\cdot)$}. We claim that $\mathcal{Z}_1$ admits an integer programming (IP) representation. This result is shown in Proposition \ref{prop:cc1}.
\begin{prop}\label{prop:cc1}
The set  $\mathcal{Z}_1$ admits an IP representation with following structure:
 \begin{equation}\label{eq:newrep}
     \mathcal{Z}_1= \left\{\bm{X} \in \mathcal{D}: x_{jt} \leq u_{jt}, \quad j \in \mathcal{J}, \; t \in \mathcal{T}\right\},
 \end{equation}
where 
\begin{align*}
    u_{jt} &=   \mathbb{I}\left [\inf_{\mathbb{P} \in \mathcal{P}_j } \mathbb{P}\biggl\{ \tilde{\Omega}: \;  [t-\tilde{\omega}_j]_{+} \leq  \rho \biggr\} \geq 1-\epsilon  \right ], \; j \in \mathcal{J}, \; t \in \mathcal{T}. 
\end{align*} 
which is a constant binary value that can be pre-computed.
\end{prop}

\begin{proof}
We first notice the following:
\begin{equation*}
    \inf_{\mathbb{P} \in \mathcal{P} } \mathbb{P}\biggl\{ \tilde{\Omega}: \;  \sum_{t \in \mathcal{T}} [t-\tilde{\omega}_j]_{+}  x_{jt}\leq  \rho \biggr\} = \inf_{\mathbb{P} \in \mathcal{P}_j } \mathbb{P}\biggl\{ \tilde{\Omega}: \;  \sum_{t \in \mathcal{T}} [t-\tilde{\omega}_j]_{+}  x_{jt}\leq  \rho \biggr\},
\end{equation*}
which holds because the summation only depends on $\tilde{\omega}_j$. For a given $\bm{X} \in \mathcal{D}$, we know that $x_{jt^*}=1$ for only one $t^* \in \mathcal{T}$. Let $\bm{X} \in \mathcal{D}$ such that $x_{jt}=1 $. It follows that, $x  \in \mathcal{Z}_1$  if and only if it satisfies the following condition:
\begin{equation*}
\inf_{\mathbb{P} \in \mathcal{P}_j } \mathbb{P}\biggl\{ \tilde{\Omega}: \;  [t-\tilde{\omega}_j]_{+} \leq  \rho \biggr\} \geq 1-\epsilon.
 \end{equation*}
We can define $u_{jt}$ as follows:
\begin{align*}
    u_{jt} &=   \mathbb{I}\left [\inf_{\mathbb{P} \in \mathcal{P}_j } \mathbb{P}\biggl\{ \tilde{\Omega}: \;  [t-\tilde{\omega}_j]_{+} \leq  \rho \biggr\} \geq 1-\epsilon  \right ], \; j \in \mathcal{J}, \; t \in \mathcal{T}. 
\end{align*} 
which is a constant value that can be pre-computed. We then notice that $x_{jt}$ belongs to $\mathcal{Z}_1$ if and only if $x_{jt} \leq  u_{jt}$, which means that \eqref{eq:newrep} provides an alternative representation of $\mathcal{Z}_1$. This completes the proof.\\
\end{proof}

\begin{rem} It is worth noting that $u_{jt}$ is non-increasing with respect to $t$. So, if $u_{j,\tau}=0$, it follows that  $u_{j t}=0,\; t\geq \tau$. Therefore, we do not need to compute $u_{j t}$ for every single $t \in \mathcal{T}$. Instead, we can use binary search to determine  $t^*_j= \min\{t: \: t \in \mathcal{T}, \; u_{j,t}=0\}$. Then, $u_{jt}= \mathbb{I}(t<t^*_j), \;j \in \mathcal{J},\;t \in \mathcal{T}$. Furthermore, when imposing $\mathcal{Z}_1$, we dictate that $x_{jt} \leq 1, \; t < t^*_j$ and that $x_{jt}=0, \; t \geq t^*_j$. Therefore, we can remove the variables $x_{jt} $ with $ t \geq t^*_j$ from the optimization problem without affecting the solution of the optimization problem.
\end{rem}

 \subsubsection{Limiting the Number of Unexpected Failures:} Unexpected failures should be avoided in industrial systems as they may induce secondary failures and can put at risk the operation of the machines. We restrict the total unexpected failures within the planning horizon by enforcing the following  DRCC set:
 \begin{equation*}
\begin{split}
        \mathcal{Z}_2 &=  \biggl\{\bm{X} \in \mathcal{D}: \inf_{\mathbb{P} \in \mathcal{P} } \mathbb{P}\biggl\{ \tilde{\Omega}: \; \sum_{j \in \mathcal{J}}\sum_{t \in \mathcal{T}} \mathbb{I}(\tilde{\omega}_j \leq t) x_{jt}   \leq \gamma \biggr\} \geq 1-\beta \biggr\}.
\end{split}
\end{equation*}
Notice that an unexpected failure occurs when a component is repaired after its failure time, i.e.,  $\tilde{\omega}_j \leq t$. Hence, the previous expression captures the total number of unexpected failures. We point out that $\mathcal{Z}_2$ can be rewritten in terms of Bernoulli random variables, resulting in the following expression:
\begin{equation}\label{eq:z2_def}
\begin{split}
        \mathcal{Z}_2 &=  \biggl\{\bm{X} \in \mathcal{D}: \inf_{\mathbb{P} \in \mathcal{P} } \mathbb{P}\biggl\{ \tilde{\Omega}: \; \sum_{j \in \mathcal{J}}\sum_{t \in \mathcal{T}} \tilde{\lambda}_{jt} x_{jt}   \leq \gamma \biggr\} \geq 1-\beta \biggr\},
\end{split}
\end{equation}
where $\tilde{\lambda}_{jt}=\mathbb{I}(\tilde{\omega}_j \leq t)$  is a Bernoulli random variable with parameter $p_{jt} = \mathbb{P}\{\tilde{\lambda}_{jt}=1\}=\mathbb{P}\{\tilde{\omega}_j \leq t\}$. This new definition of $\mathcal{Z}_2$ is useful to derive some structural properties of this DRCC set. To this end, we show two preliminary Lemmas next.

\begin{lemma}\label{lemma:1}
Let $(x_i,y_i), \; i=1,...,N,$ be pairs of independent Bernoulli random variables with parameters $0\leq p_i^x \leq 1$ and $0 \leq p_i^y \leq 1$, respectively. If $p_i^x \leq p_i^y, \; i=1,...,N,$ then  
 \begin{equation}\label{eq:order}
    \mathbb{P} \left \{\sum_{i=1}^{N}x_i \leq  u \right\} \geq \mathbb{P}\left\{\sum_{i=1}^{N}y_i \leq  u \right\}, \; u  \in \mathbb{Z}.
\end{equation}
\end{lemma}
\begin{proof}
When $p_i^x \leq p_i^y, \; i=1,...,N$, \cite{kochar2022stochastic}  showed that $\sum_{i=1}^{N}x_i \leq_{lr} \sum_{i=1}^{N}y_i $, where $\leq_{lr}$ represents the \textit{likelihood ratio order} between these two random variables. From the fact that $\leq_{lr} \Rightarrow \leq_{st}$, it follows that $\sum_{i=1}^{N}x_i \leq_{st} \sum_{i=1}^{N}y_i$, where $\leq_{lr}$ represents the \textit{usual stochastic order} between these two random variables. From the definition of $\leq_{lr}$, we conclude that $    \mathbb{P} \left\{\sum_{i=1}^{N}x_i > u \right\} \leq \mathbb{P}\left\{\sum_{i=1}^{N}y_i > u \right\}.$ Finally, we obtain the desired result \eqref{eq:order} using the identity $ \mathbb{P} \left \{  > u \right\}=1-\mathbb{P} \left(  \leq  u \right)$.
\end{proof}

 \begin{lemma}\label{lemma:bernoulli_cdf}
Let $x_1,x_2,...,x_n$ be independent Bernoulli random variables with parameter $p_i, \; i=1,...,n$, respectively. Let $m$ be a positive integer. Then, $\mathbb{P}\left\{ \sum_{i=1}^n x_i \leq m \right\} = v_{nm}$, where $v_{nm}$ can be computed by solving the following recursive equations defined in terms of $v_{jk} \in [0,1], \; j=0,1,...,n,\; k=0,1,...,m$.
\begin{subequations}
  \begin{empheq}[left={v_{jk}=\empheqlbrace}]{align}
        &1 & \text{if }j \leq  k\geq 0 \quad \quad \quad \quad \; \; \label{eq:rec1} \\
        &v_{j-1,0}  (1-p_j) & \text{elseif } j >   k=0. \quad \quad \quad \label{eq:rec2}\\
        &v_{j-1,k-1}  p_j +v_{j-1,k}  (1-p_j)  & \text{elseif }j> k>0 \quad \quad \quad \label{eq:rec3} 
  \end{empheq}
  \end{subequations}
 \end{lemma}
 \begin{proof} Let us define $v_{jk}  \; j=0,1,...,n,\; k=0,1,...,m$ as follows:
 \begin{equation}\label{eq:v_def}
     v_{jk}=\mathbb{P}\left\{ \sum_{i=1}^j x_i \leq k \right\}.
 \end{equation}
 From its definition, it is direct that $\mathbb{P}\left\{ \sum_{i=1}^n x_i \leq m \right\} = v_{nm}$. So, it only remains to show that $v_{nm}$ can be computed using  \eqref{eq:rec1}-\eqref{eq:rec3}. For this, we will show that $v_{jk}$ defined in \eqref{eq:v_def} satisfies \eqref{eq:rec1}-\eqref{eq:rec3} and that the system of equations \eqref{eq:rec1}-\eqref{eq:rec3} has an unique solution.
 
 First, if $j\leq k \geq 0$,  $v_{jk}$ satisfies \eqref{eq:rec1}. This follows from the fact that $x_i$ is binary, implying that $\sum_{i=1}^j x_i \leq j$ with probability 1. Second, if  $j> k=0$, $v_{j0}$ can be rewritten as follows:
\begin{align*}
v_{j0}&=\mathbb{P}\left\{ \sum_{i=1}^j x_i \leq 0 \right\}\\
    & =  \mathbb{P}\left\{ \sum_{i=1}^{j-1} x_i \leq 0 | x_j=0\right\} \mathbb{P}\left\{ x_j=0 \right\}\\
    & =  \mathbb{P}\left\{ \sum_{i=1}^{j-1} x_i \leq 0  \right\} \mathbb{P}\left\{ x_j=0 \right\}\\
    &= v_{j-1,0} (1-p_j),
\end{align*}
 thus $v_{j0}$ satisfies \eqref{eq:rec2}. Third, if $j>k>0$ $v_{j0}$ can be rewritten as follows:
\begin{align*}
v_{jk}&=\mathbb{P}\left\{ \sum_{i=1}^j x_i \leq k \right\}\\
    & = \mathbb{P}\left\{ \sum_{i=1}^{j-1} x_i \leq k-1 | x_j=1\right\} \mathbb{P}\left\{ x_j=1 \right\} +\mathbb{P}\left\{ \sum_{i=1}^{j-1} x_i \leq k | x_j=0\right\} \mathbb{P}\left\{ x_j=0 \right\}\\
    & = \mathbb{P}\left\{ \sum_{i=1}^{j-1} x_i \leq k-1\right\} p_j +\mathbb{P}\left\{ \sum_{i=1}^{j-1} x_i \leq k  \right\} (1-p_j)\\
    &= v_{j-1,k-1} p_j+ v_{j-1,k} (1-p_j),
\end{align*}
so $v_{jk}$ satisfies \eqref{eq:rec3}.

Finally, we notice that the equations \eqref{eq:rec1}-\eqref{eq:rec3} are recursive and have a unique solution once their initial condition is fixed. As we are fixing the initial condition with \eqref{eq:rec1}, we can conclude that $v_{jk}$ defined in \eqref{eq:v_def} is the unique solution to this system of equations. This finishes the proof.\\
\end{proof}
 \begin{eg}
Let $x_1,x_2,x_3$ be independent Bernoulli random variables with parameter $p_1=0.1,p_2=0.25, p_3=0.8$, respectively. Let us compute $\mathbb{P}\left\{ \sum_{i=1}^n x_i \leq 2 \right\}$ using Lemma \ref{lemma:bernoulli_cdf}. First, from \eqref{eq:rec1} it follows that $v_{00}=v_{01}=v_{11}=v_{02}=v_{12}=v_{22}=1$. Second, from \eqref{eq:rec1} we obtain that $v_{10}=v_{00}(1-p_1)=0.9$, $v_{20}=v_{10}(1-p_2)=0.675$, $v_{30}=v_{20}(1-p_3)=0.135$. Third, from \eqref{eq:rec1} we obtain that  $v_{21}=v_{10}p_2+v_{11}(1-p_2)=0.975$, $v_{31}=v_{20}p_3+v_{21}(1-p3)=0.735$, and $v_{32}=v_{21}p_3+v_{22}(1-p_3)=0.98.$
 \end{eg}

\begin{prop}\label{prop:cc2}
 The DRCC set $\mathcal{Z}_2$ can be rewritten as a chance constrained set, i.e.
\begin{equation}\label{eq:z_2_cc}
\begin{split}
        \mathcal{Z}_2 &=  \biggl\{\bm{X} \in \mathcal{D}:   \mathbb{P}\biggl\{ \tilde{\Omega}: \; \sum_{j \in \mathcal{J}}\sum_{t \in \mathcal{T}} \Tilde{\eta}_{jt} x_{jt}   \leq \gamma \biggr\} \geq 1-\beta \biggr\},
\end{split}
\end{equation}
where $\Tilde{\eta}_{jt}$  is a Bernoulli random variable with $\bar{p}_{jt} = \mathbb{P}(\Tilde{\eta}_{jt}=1)= \sup_{\mathbb{P} \in \mathcal{P}_j } \mathbb{P}(\tilde{\omega}_j \leq t), \; j \in \mathcal{J}, \; t \in \mathcal{T}.$
\end{prop}
\begin{proof}
It suffices to show that
\begin{equation}\label{eq:part0}
      \inf_{\mathbb{P} \in \mathcal{P} } \mathbb{P}\biggl\{\tilde{\Omega}: \; \sum_{j \in \mathcal{J}}\sum_{t \in \mathcal{T}} \tilde{\lambda}_{jt} x_{jt}   \leq \gamma \biggr\} =\mathbb{P}\biggl\{\tilde{\Omega}: \; \sum_{j \in \mathcal{J}}\sum_{t \in \mathcal{T}} \Tilde{\eta}_{jt} x_{jt}   \leq \gamma \biggr\},\; \bm{X} \in \mathcal{D}.
\end{equation}
First, from the definition of $\inf_{\mathbb{P} \in \mathcal{P} }(\cdot)$ it follows that  
\begin{equation}\label{eq:part1}
\inf_{\mathbb{P} \in \mathcal{P} } \mathbb{P}\biggl\{\tilde{\Omega}: \; \sum_{j \in \mathcal{J}}\sum_{t \in \mathcal{T}} \tilde{\lambda}_{jt} x_{jt}   \leq \gamma \biggr\} \leq \mathbb{P}\biggl\{\tilde{\Omega}: \; \sum_{j \in \mathcal{J}}\sum_{t \in \mathcal{T}} \Tilde{\eta}_{jt} x_{jt}   \leq \gamma \biggr\}  ,\; \bm{X} \in \mathcal{D}.
\end{equation}
Secondly, we pick arbitrary $\mathbb{P}' \in \mathcal{P} $ and $\bm{X} \in \mathcal{D}$. As $\bm{X} \in \mathcal{D}$, it implies that for every $j\in \mathcal{J}$ $x_{jt}=1$ for a unique $t \in \mathcal{T}$. Let us define $\tau_j \in \mathcal{T}$ such that $x_{j \tau_j}=1, \; j \in \mathcal{J}$. Using this definition, we obtain 
\begin{subequations}
\begin{align}
\mathbb{P}\biggl\{\tilde{\Omega}: \; \sum_{j \in \mathcal{J}}\sum_{t \in \mathcal{T}} \tilde{\lambda}_{jt} x_{jt}   \leq \gamma \biggr\} &=\mathbb{P}\biggl\{\tilde{\Omega}: \; \sum_{j \in \mathcal{J}}  \tilde{\lambda}_{j, \tau_j}    \leq \gamma \biggr\} \label{eq:part2_1}\\
& \geq  \mathbb{P}\biggl\{\tilde{\Omega}: \; \sum_{j \in \mathcal{J}}  \Tilde{\eta}_{j, \tau_j}    \leq \gamma \biggr\}   \label{eq:part2_2}\\
&= \mathbb{P}\biggl\{\tilde{\Omega}: \; \sum_{j \in \mathcal{J}}\sum_{t \in \mathcal{T}} \Tilde{\eta}_{j t}  x_{jt}   \leq \gamma \biggr\},\label{eq:part2_3}
\end{align}
\end{subequations}
where the first inequality is because $\sum_{j \in \mathcal{J}}  \eta_{j, \tau_j}$ is a sum of independent Bernoulli random variables with parameter $p_{j,\tau_j}= \mathbb{P}'\{\tilde{\omega}_j \leq \tau_j\} \leq \bar{p}_{j\tau_j} = \sup_{\mathbb{P} \in \mathcal{P} } \mathbb{P}\{\tilde{\omega}_j \leq \tau_j\}$ and thus the result in Lemma \ref{lemma:bernoulli_cdf} follows, and the second equality is obtained from the definition of $\tau_j$.
As \eqref{eq:part2_3} holds any $\mathbb{P}' \in \mathcal{P} $ and $\bm{X} \in \mathcal{D}$, it follows that for any $\bm{X}\in \mathcal{D}$
\begin{equation}\label{eq:part3}
\inf_{\mathbb{P} \in \mathcal{P} }  \mathbb{P}\biggl\{\tilde{\Omega}: \; \sum_{j \in \mathcal{J}}\sum_{t \in \mathcal{T}} \tilde{\lambda}_{jt} x_{jt}   \leq \gamma \biggr\}  \geq \mathbb{P}\biggl\{\tilde{\Omega}: \; \sum_{j \in \mathcal{J}}\sum_{t \in \mathcal{T}} \Tilde{\eta}_{j t}  x_{jt}   \leq \gamma \biggr\}.
\end{equation}
Combining the results of \eqref{eq:part1} and \eqref{eq:part3}, we obtain the equality of \eqref{eq:part0}, thus finishing the proof.
\end{proof}

\begin{prop}\label{prop:refor_cc2}
    The DRCC set $\mathcal{Z}_2$ admits the following MILP representation:
{\small\begin{equation}\label{eq:final_Z2}  \hspace*{-.5cm}   \mathcal{Z}_2=
    \left\{
    \begin{array}{ll}
    \begin{array}{l}
         \bm{X} \in \mathcal{D}\\ 
         \bm{a} \in [0,1]^{J\times \gamma}\\
         \bm{b} \in [0,1]^{J\times T} \\
         \bm{c} \in [-1,0]^{J\times \gamma \times T}
      \end{array} \vline
      \begin{array}{ll}
      a_{J,\gamma}\geq 1-\beta   \\[5pt]
       a_{je}=\left \{
      \begin{array}{lll}
        1     \\
         a_{j-1,0}- \sum_{t \in \mathcal{T}} b_{jt}   \\
        a_{j-1,e}  + \sum_{t \in \mathcal{T}} c_{jet}      
      \end{array} \right. &  \begin{split}
     &\text{if }j \leq  e\geq 0  \\
      &    \text{elseif } j >   e=0 \\
      & \text{elseif }j> e>0  
      \end{split} \\[20pt]
    0  \leq b_{jt}   \leq  \bar{p}_{jt} x_{jt} & \text{if } j >   e=0  \\
    0  \leq a_{j-1,0} \bar{p}_{jt}- b_{jt} \leq \bar{p}_{jt}(1-x_{jt}) & \text{if } j >   e=0  \\[10pt]
    -\bar{p}_{jt} x_{jt}   \leq c_{jet}   \leq  0 & \text{if }j> e>0 \\
   -\bar{p}_{jt}(1-x_{jt}) \leq \bar{p}_{jt}   (a_{j-1,e-1}-a_{j-1,e})  - c_{jet} \leq 0 & \text{if }j> e>0 
      \end{array}
    \end{array}
    \right\}.
\end{equation}}
\end{prop}

\begin{proof}
Using Proposition \ref{prop:cc2} we know that $\mathcal{Z}_2$ admits the following representation:
\begin{equation*} 
\begin{split}
        \mathcal{Z}_2 &=  \biggl\{\bm{X} \in \mathcal{D}:   \mathbb{P}\biggl\{\tilde{\Omega}: \; \sum_{j \in \mathcal{J}}\sum_{t \in \mathcal{T}} \Tilde{\eta}_{jt} x_{jt}   \leq \gamma \biggr\} \geq 1-\beta \biggr\}.
\end{split}
\end{equation*}
For a given $\bm{X} \in \mathcal{D}$, it verifies that for every $j\in \mathcal{J}$ $x_{jt}=1$ for an unique $t \in \mathcal{T}$. Let us define $\tau_j \in \mathcal{T}$ such that $x_{j \tau_j}=1$. Thus,

\begin{equation*}
    \mathbb{P}\biggl\{\tilde{\Omega}: \; \sum_{j \in \mathcal{J}}\sum_{t \in \mathcal{T}} \Tilde{\eta}_{jt} x_{jt}   \leq \gamma \biggr\}   =   \mathbb{P}\biggl\{\tilde{\Omega}: \;  \sum_{j \in \mathcal{J}} \Tilde{\eta}_{j\tau_j}  \leq \gamma \biggr\}  
\end{equation*}
where  $ \sum_{j \in \mathcal{J}} \Tilde{\eta}_{j\tau_j}$  is a sum of $J$ independent Bernoulli random variables with parameters $p_j=\sum_{t\in\mathcal{T}}\bar{p}_{jt} x_{jt}, \; j\in \mathcal{J}$. Therefore, we can compute $  \mathbb{P}\biggl\{\tilde{\Omega}: \;  \sum_{j \in \mathcal{J}} \Tilde{\eta}_{j\tau_j}  \leq \gamma \biggr\}  $ using Lemma \ref{lemma:bernoulli_cdf}. For this we introduce new decision variables $a_{je} \in [0,1], \; j=0,1,...,J,\; e=0,1,...,\gamma$, which must satisfy the following constraint:
 
  \begin{empheq}[left={a_{je}=\empheqlbrace}]{align}
        &1 & \text{if }j \leq  e\geq 0 \quad \quad \quad \quad \; \;   \nonumber\\
        &a_{j-1,0}-a_{j-1,0}  p_j & \text{elseif } j >   e=0 \quad \quad \quad  \nonumber\\
        &(a_{j-1,e-1}-a_{j-1,e})  p_j + a_{j-1,e}  & \text{elseif }j> e>0 \quad \quad \quad \nonumber 
  \end{empheq} 
If $j >   e=0$, we introduce $b_{jt} \in \mathbb{R}, t\in \mathcal{T}$. Then, 
\begin{align*}
    a_{je} &= a_{j-1,0}-a_{j-1,0}  p_j \\
     & = a_{j-1,0}- \sum_{t \in \mathcal{T}} \bar{p}_{jt} x_{jt}  a_{j-1,0}\\
      & = a_{j-1,0}- \sum_{t \in \mathcal{T}} b_{jt},
\end{align*} 
where $b_{jt}$ satisfies the following constraints:
\begin{equation*}
\begin{split}
    0 & \leq b_{jt}   \leq  \bar{p}_{jt} x_{jt} \\
    0 &\leq a_{j-1,0} \bar{p}_{jt}- b_{jt} \leq \bar{p}_{jt}(1-x_{jt})
\end{split}, \; j \in \mathcal{J}, \; t \in \mathcal{T}.
\end{equation*}
If $j> e>0$, we introduce $c_{jet}\in \mathbb{R}, \; t \in \mathcal{T}$
 
\begin{align*}
    a_{je} &= (a_{j-1,e-1}-a_{j-1,e})  p_j + a_{j-1,e}  \\
     & = a_{j-1,e} + \sum_{t \in \mathcal{T}} \bar{p}_{jt} x_{jt}   (a_{j-1,e-1}-a_{j-1,e}) \\
      & = a_{j-1,e}  + \sum_{t \in \mathcal{T}} c_{jet},
\end{align*} 
where $c_{jet}$ satisfies the following constraints:
\begin{equation*}
\begin{split}
    -\bar{p}_{jt} x_{jt} & \leq c_{jet}   \leq  0 \\
   -\bar{p}_{jt}(1-x_{jt})&\leq \bar{p}_{jt}   (a_{j-1,e-1}-a_{j-1,e})  - c_{jet} \leq 0
\end{split}, \; j \in \mathcal{J}, \; j> e>0, \; t \in \mathcal{T}.
\end{equation*}
Combining the previous steps, we end up with the MILP reformulation presented in \eqref{eq:final_Z2}, thus finishing the proof.
\end{proof}

 \begin{cor}\label{cor:reformulation}
The proposed DRCC program presented in \eqref{eq:objective} can be reformulated exactly as a MILP.  
\end{cor}
\begin{proof}
Direct from   Propositions 1, 2, and 4.\\
\end{proof}
\edit{The MILP reformulation obtained in Corollary \ref{cor:reformulation} relies on the ambiguity set $\mathcal{P}$ defined in \eqref{eq:global_ambiguity_set}, which is a generic discrepancy-based ambiguity set. Thus, the result of Corollary \ref{cor:reformulation} is quite general. This reformulation can be easily implemented in off-the-shelf solvers, as long as the reformulation parameters ($\psi_{jt}, u_{jt}, \bar{p}_{jt},\ j \in \mathcal{J}, \; t \in \mathcal{T}$) can be precomputed efficiently. However, the computation of these parameters depends on the specific ambiguity set adopted.
}

\subsection{Derivations when adopting a type-$\infty$ Wasserstein ambiguity set}\label{subsec:wasserstein_reformulation}

\edit{In this section, we examine a specific ambiguity set to demonstrate how the parameters, $\psi_{jt}, u_{jt}, \bar{p}_{jt},\ j \in \mathcal{J}, \; t \in \mathcal{T}$, can be determined in practice. To this end, we adopt a type-$\infty$ Wasserstein ambiguity set, primarily for its mathematical tractability. This approach enables us to derive exact formulas for efficiently computing the parameters of the MILP reformulation as shown in Proposition \ref{prop:infty_reformulation}. }

We first construct an empirical distribution of the remaining lifetime $\tilde{\omega}_j$  by generating $N$ independent remaining lifetime samples from the distribution $\widehat{\mathrm{RLD}}_j$. The set of samples is indexed by $ i \in  \mathcal{N}=\{1,...,N\}$. The $i$-th sample is denoted as $\omega_j^i$. We define the empirical RLD for component $j$ as follows:
\begin{equation*}
\widehat{\mathbb{P}}_N^j (\tilde{\omega}_j)= \frac{1}{N}\sum_{i \in \mathcal{N}} \mathbb{I}(\tilde{\omega}_j=\omega_j^i).
\end{equation*}
Then, $\infty-$Wasserstein ambiguity set \citep{xie2021distributionally} is defined as follows: 
\begin{equation}\label{eq:infinity_ambiguity}
   \mathcal{P}_j =  \left \{ \mathbb{P}: \; \mathbb{P}\left \{\tilde{\omega} \in \Xi_j \right\}=1,\; W_{\infty} \left(\mathbb{P}, \widehat{\mathbb{P}}_N^j (\tilde{\omega}_j) \right) \leq \delta_j \right\},
\end{equation}
where  the $\infty-$Wasserstein distance is defined as
\begin{equation*}  
\begin{split}
W_{\infty} \left(\mathbb{P}_1, \mathbb{P}_2\right) = &  \inf_{\mathrm{Q}} \biggl\{
 \text{ess.sup} \Vert \tilde{\omega}_1-\tilde{\omega}_2 \Vert  : \begin{array}{l}
\text{$\mathrm{Q}$ is a joint distribution of  $\tilde{\omega}_1$ and $\tilde{\omega}_2$} \\
\text{with marginals $\mathbb{P}_1$ and $\mathbb{P}_2$, respectively}
\end{array}\biggr\}.
\end{split}
\end{equation*}
Therefore, in this case, $\mathcal{P}_j$ corresponds to the set of all distributions taking values in $\Xi_j$ whose $\infty-$Wasserstein distance to the empirical distribution $\widehat{\mathbb{P}}_N^j$ is less or equal than $\delta_j$. This set constitutes a Wasserstein ball with radius $\delta_j$ and center $\widehat{\mathbb{P}}_N^j$.  

\begin{prop}\label{prop:infty_reformulation}
If $\mathcal{P}_j, \; j \in \mathcal{J}$ is defined as a type-$\infty$ Wasserstein ambiguity set \eqref{eq:infinity_ambiguity}, then:\\
\begin{enumerate}[(a)]
\item the parameters $\psi_{jt},\;  j \in \mathcal{J},\;  t \in \mathcal{T}$ used in the reformulation resulting from Proposition \ref{prop:obj} can be computed as $\psi_{jt}= \frac{1}{ N } \sum_{i \in \mathcal{N}}\theta_{ijt}$, where
\begin{align*}
 \theta_{ijt} =  \max\Biggl\{ \left(C_j^{pr}+V_j^{pr}(\min \left\{\omega_j^i+ \delta_j , \max{\Xi_j}  \right\} -t) \right)
 \mathbb{I}\left (t<\omega_j^i+\delta_j    \right), \\  \left(C_j^{co}+V_j^{co}(t-\max \{\omega_j^i- \delta_j , 0 \}) \right) \mathbb{I}\left ( \omega_j^i-\delta_j \leq t  \right)  \Biggr\}.
\end{align*}

\item the parameters $u_{jt},\;  j \in \mathcal{J},\;  t \in \mathcal{T}$ used in the reformulation resulting from Proposition \ref{prop:cc1} can be computed as
\begin{equation*}
 u_{jt} = \mathbb{I}\left[ \frac{1}{ N }\sum_{i \in \mathcal{N}} \mathbb{I}\left\{ [t-\max(\omega_j^i-\delta ,0)]_{+}\leq  \rho \right\} \geq 1-\epsilon  \right] , \quad  j \in \mathcal{J}, \; t \in \mathcal{T}.
\end{equation*}

\item the parameters $\bar{p}_{jt},\;  j \in \mathcal{J},\;  t \in \mathcal{T}$ used in the reformulation resulting from Proposition \ref{prop:cc2} can be computed as
\begin{equation*}
 \bar{p}_{jt}  =\frac{1}{N} \sum_{i\in \mathcal{N}} \mathbb{I} \left(\omega_j^i- \delta_j  \leq t \right) , \quad  j \in \mathcal{J}, \; t \in \mathcal{T}.
\end{equation*}
\end{enumerate}
    
\end{prop}

\begin{proof}
First, we show (a). We use the result of Proposition 3 in \citep{bertsimas2022data} to derive the following equality:
\begin{equation}\label{eq:alpha_inf}
\sup_{\mathbb{P} \in \mathcal{P}_j}{\mathbb{E}_{\mathbb{P}}   \left[\alpha(\tilde{\omega}_j,t) \right]}  = \frac{1}{ N }  \sum_{i\in \mathcal{N}} \sup_{
\omega_j: |\omega_j - \omega_j^i|\leq \delta_j}  \alpha(\omega_j,t).    
\end{equation} 
Subsequently, we define $\theta_{ijt}$ as follows: 
\begin{equation}\label{eq:solve_alpha1}
\begin{split}
       &\theta_{ijt}=\sup_{\omega_j: | \omega_j- \omega_j^i| \leq \delta_j }  \alpha(\omega_j,t)  = \max \left\{   \begin{split}
\max \; &  \; V^{pr}_j(\omega_j-t)  + C^{pr}_j \\
\text{st. } & \max{\Xi_j} \geq \omega_j > t \\
& | \omega_j- \omega_j^i| \leq \delta_j 
\end{split} , \quad \begin{split}
\max \; &  \;  V^{co}_j(t-\omega_j)+ C^{co}_j \\
\text{st. } & 0 \leq \omega_j \leq   t \\
& | \omega_j- \omega_j^i| \leq \delta_j 
\end{split} \right\}.
\end{split}
\end{equation}
The right-hand side of \eqref{eq:solve_alpha1} can be further simplified as follows:
\begin{align*}
 \theta_{ijt} =  \max\Biggl\{ \left(C_j^{pr}+V_j^{pr}(\min \left\{\omega_j^i+ \delta_j , \max{\Xi_j}  \right\} -t) \right)
 \mathbb{I}\left (t<\omega_j^i+\delta_j     \right), \\  \left(C_j^{co}+V_j^{co}(t-\max \{\omega_j^i- \delta_j , 0 \}) \right) \mathbb{I}\left ( \omega_j^i-\delta_j \leq t  \right)  \Biggr\}.
\end{align*}
The proof of (a) is finished by taking: $\psi_{jt}= \frac{1}{ N } \sum_{i \in \mathcal{N}}\theta_{ijt}$.

Second, we show (b). \edit{We derive \eqref{eq:bb1} using the classical equality between probability and the expected value of the indicator function.} We then obtain \eqref{eq:bb2} using the result of Proposition 3 from \citep{bertsimas2022data}. Next, we use the fact that the indicator function is monotonous in $\omega_j$ to derive  \eqref{eq:bb3}, where the supremum can be easily calculated leading to \eqref{eq:bb4} and concluding the proof.
\begin{subequations}
\begin{align}
    u_{jt} = \mathbb{I}\left [\inf_{\mathbb{P} \in \mathcal{P}_j } \mathbb{P}\biggl\{   [t-\tilde{\omega}_j]_{+} \leq  \rho \biggr\} \geq 1-\epsilon  \right ] &= \mathbb{I}\left [\inf_{\mathbb{P} \in \mathcal{P}_j } \mathbb{E}^{\mathbb{P}}\left[\mathbb{I}\biggl\{   [t-\tilde{\omega}_j]_{+} \leq  \rho \biggr\} \right] \geq 1-\epsilon  \right ] \label{eq:bb1}\\
    & \mkern-18mu=   \mathbb{I}\left [ \frac{1}{ N }\sum_{i \in \mathcal{N}} \inf_{\omega_j: |\omega_j-\omega_j^i| \leq \delta_j } \mathbb{I}\left\{ [t-\omega_j]_{+}\leq  \rho \right\} \geq 1-\epsilon  \right ] \label{eq:bb2}\\ 
    &\mkern-18mu=   \mathbb{I}\left [ \frac{1}{ N }\sum_{i \in \mathcal{N}}  \mathbb{I}\left\{ \sup_{\omega_j: |\omega_j-\omega_j^i| \leq \delta_j }[t-\omega_j]_{+}\leq  \rho \right\} \geq 1-\epsilon  \right ] \label{eq:bb3}\\ 
    &\mkern-18mu= \mathbb{I}\left[ \frac{1}{ N }\sum_{i \in \mathcal{N}} \mathbb{I}\left\{ [t-\max(\omega_j^i-\delta ,0)]_{+}\leq  \rho \right\} \geq 1-\epsilon  \right]. \label{eq:bb4}
\end{align}
\end{subequations}

Third, we show (c). We obtain \eqref{eq:identity} using the result of Proposition 3 from \citep{bertsimas2022data}. Next, we use the fact that the indicator function is monotonous in $\omega_j$ to derive  \eqref{eq:final_result1}, where the infimum can be easily computed leading to \eqref{eq:final_result2} and concluding the proof.
\begin{subequations}
\begin{align}
 \bar{p}_{jt} = \sup_{\mathbb{P} \in \mathcal{P}_j} \mathbb{P}\{\tilde{\omega}_j \leq t\}  =  \sup_{\mathbb{P} \in \mathcal{P}_j} \mathbb{E}_{\mathbb{P}} \left[ \mathbb{I}\left( \tilde{\omega}_j\leq t \right) \right]& =\frac{1}{ N }  \sum_{i\in \mathcal{N}}  \sup_{ \omega_j: |  \omega_j-  \omega_j^i|\leq \delta_j } \mathbb{I}\left( \omega_j\leq t \right) \label{eq:identity}\\
 & \mkern-10mu = \frac{1}{N} \sum_{i\in \mathcal{N}} \mathbb{I} \left( \inf_{ \omega_j: |  \omega_j-  \omega_j^i|\leq \delta_j } \omega_j  \leq t \right) \label{eq:final_result1}\\
  & \mkern-10mu = \frac{1}{N} \sum_{i\in \mathcal{N}} \mathbb{I} \left(\omega_j^i- \delta_j  \leq t \right). \label{eq:final_result2}
\end{align}
\end{subequations}
\end{proof}
Therefore, when adopting a type-$\infty$ Wasserstein ambiguity set, Proposition \ref{prop:infty_reformulation} provides exact formulas for efficiently computing the parameters of the MILP reformulation. \revise{While the type-$\infty$ Wasserstein ambiguity is convenient from a mathematical tractability perspective, its interpretation may not be straightforward from its definition in \eqref{eq:infinity_ambiguity}. However, for bounded real-valued random variables, the type-$\infty$ Wasserstein distance can be written as \citep{ramdas2017wasserstein}:
$$W_\infty(\mathbb{P},\mathbb{Q}) = \Vert F^{-1} - G^{-1} \Vert_\infty,$$
where $F^{-1}$ and $G^{-1}$ are the quantile functions of $\mathbb{P}$ and $\mathbb{Q}$, and $\Vert \cdot \Vert_\infty$ is the supremum norm. Thus, $W_\infty$ captures the maximum difference between quantile functions. This formulation applies to our setting since the ambiguity sets $\mathcal{P}_j,\; j \in \mathcal{J}$ in \eqref{eq:infinity_ambiguity} involve bounded real-valued random variables.

Our numerical studies utilize the type-$\infty$ Wasserstein ambiguity, leveraging the tractable reformulation derived in this section. However, it is important to note that adopting the type-$\infty$ Wasserstein distance may introduce certain limitations. Specifically, the type-$p$ Wasserstein distance satisfies $W_p(\mathbb{P}, \mathbb{Q}) \leq W_\infty(\mathbb{P}, \mathbb{Q})$ for all $p \geq 1$ \citep{santambrogio2015optimal}. Therefore, for the same radius $\delta$, a DRO formulation based on type-$\infty$ ambiguity may yield more conservative solutions than one based on type-$p$ Wasserstein sets. Future work will explore alternative ambiguity sets, including those defined via type-$p$ Wasserstein distances, moments, or $\phi$-divergences, to assess their respective advantages and limitations in solving joint optimization problems for maintenance and spare parts provisioning. }

\begin{rem} It is important to note that the proposed DRCC formulation for jointly optimizing maintenance and spare parts inventory, presented in Section \ref{sec:formulation}, can also be applied to purely maintenance scheduling. To do this, it suffices to remove the decision variables ($h_{lt}, g_{lt}^{reg}, g_{lt}^{exp}, r_t, \; l \in \mathcal{L}, \; t \in \mathcal{T}$) and constraints (\eqref{eq:fixed_order} and \eqref{eq:inventory}) related to spare parts inventory. The resulting CBM model is simpler, with fewer decision variables and constraints, yet it retains the DR approach and admits the same MILP reformulation presented in Section \ref{sec:formulation}.

\end{rem}

\section{Numerical Studies}\label{sec:computational_studies}

Our numerical study is inspired by a wind farm application \citep{bakir2021integrated}. We rely on a set of simulated CM data to represent component degradation. Characteristics of the simulated data were inspired by a real-world vibration monitoring dataset used in \citep{gebraeel2005residual}. We assume that we have a wind farm comprised of 50 wind turbines. Each turbine comprises three different types of critical components. \edit{ The expected lifetimes of these three components are 6, 8, and 9 months}. For example, these components can be a gearbox, converter, and main turbine bearing. Thus, we expect to manage three different spare parts inventories \add{one for each type of critical component}, i.e., $\mathcal{L}=\{1,2,3\}$. \add{Our objective is to}  solve the joint optimization problem for the set of 150 components that are distributed across the 50 wind turbines, i.e., $\mathcal{K}=\{1,...,50\}$ and $\mathcal{J}=\{1,...,150\}$.

We implemented a simulation framework to emulate CM data collected from the different types of components. This simulation framework utilizes a stochastic model similar to the one presented in \citep{elwany2009real} and is described in  Appendix A. This simulation framework was used to simulate CM data from an "as good as new" state up to a failure state, where failure was defined as the trend in the CM data crossing a predefined failure threshold. The resulting $n-$th degradation signal of component type $l$ corresponds to collection of CM observations and is denoted as $\mathcal{S}_{l}^n=\{S(0),...,S(t_{ln})\}$, where $S(t)$ represents CM a observation collected at time $t$.

Given the model used to simulate the CM data, the RLD of a component type $l$ at age $t$ can be estimated using an Inverse Gaussian distribution $IG(\nu_{lt},\gamma_{lt})$. The location parameter $\nu_{lt}>0$ and the shape parameter $\gamma_{lt}>0$ can be inferred from CM data by the Bayesian updating framework presented by \cite{elwany2009real}. This updating process requires 1) a set of historical degradation signals $\mathcal{S}_l^{\text{train}}= \{\mathcal{S}_{l}^1,\mathcal{S}_{l}^2,...\}$ hereinafter referred to as \textit{training data}  used to characterize the prior distributions of $\nu_{lt}$ and $\gamma_{lt}$ and 2) a set of CM observations $\{S(0),...,S(t_{t})\}$ collected from the component under study. Consequently, estimates of $\nu_{lt}$ and $\gamma_{lt}$ can be refreshed as new observations become available.

The periodic updating process of the prognostic model parameters often helps in improving the prediction accuracy and thus motivates the implementation of the optimization model in a rolling horizon fashion. We define a freeze period $[1, \ldots, \Delta^{upd}] \subset \mathcal{T}$, where $\Delta^{upd}$ denotes the length of the freeze period.   \add{Only the portion of the decisions that fall within the freeze period is actually executed. Decisions outside the freeze period can be revised based on new CM data observations on the next resolve epoch.}  This means that only maintenance decisions scheduled within the freeze period are performed. After the freeze period, components' RLDs are re-estimated using fresh sensor data, and then the optimization model is re-executed, leading to updated decisions. This decision-updating process is repeated over time so that the optimization model can adapt its decisions to the degradation process experienced by the components.

The optimization model was implemented with following parameters: $T_{max}=50$ [days], $\Delta^{upd}=20$ [days], $G_t=30$ [units], $V_j^{pr}= \$0.13K$, $V_j^{co}= 6 V_j^{pr}$, $C_j^{pr}=\$ 1.3 K$, \edit{ $C_j^{co}= 6 C_j^{pr}$}, $C_k^{down}= \$2K$, $C^{crew}= \$40K$, $C_l^{reg}=\$1.5K$, \edit{$C_l^{exp}=4 C_l^{reg}$}, $C_l^{hol}=\$0.1K$, $B^{reg}=\$3K$, $\Delta^{reg}=20$ [days], $\rho=5$ $\gamma=7$, $\epsilon=0.1$ and $\beta=0.1$. \edit{Notice that the planning horizon is 50 days, while the expected lifetimes of the three types of components under study are 6, 8, and 9 months. Therefore, Assumptions \ref{as:one_failure} and \ref{as:one_repair} are satisfied in this setting. } We adopted a type-$\infty$ Wasserstein ambiguity set and so used the expressions derived in Section \ref{subsec:wasserstein_reformulation} to reformulate the corresponding optimization model as a MILP. For testing purposes, we defined the Wasserstein radius as $\delta_j=\delta \cdot \sigma_j, \; j \in \mathcal{J}$, where $\delta$ is referred to as the \textit{normalized radius} and  $\sigma_j$ corresponds to the standard deviation of $\tilde{\omega}_j$. We conducted 3 different studies. Section \ref{subsec:varying_traindata} evaluates the performance of the proposed DRCC formulation in two different situations considering sparse and abundant \textit{training data} for the prognostic algorithm. Section \ref{subsec:sequential_joint} compares the performance of a DRCC joint optimization for a maintenance schedule and spare provisioning versus a DRCC sequential optimization, which first optimizes the maintenance schedule and then spares provisioning. Section \ref{subsec:running_times} studies the running time of the DRCC model for two different solution methods.

\subsection{Effects of the Amount of Training Data used by the Prognostic Model}\label{subsec:varying_traindata}

We conducted two simulation studies to analyze the performance of the optimization model depending on the amount of \textit{training data} for the prognostic algorithm, i.e., depending on $|\mathcal{S}_l^{\text{train}}|$. Each simulation study included 20 rolling horizon updates, equivalent to a simulation period of 400 days. The simulations were repeated 20 times for each configuration to account for the fluctuations in the simulated environment. Specifically, we analyzed two simulation setups referred to as  \textit{Sparse training data} ( $|\mathcal{S}_l^{\text{train}}|=5$) and \textit{Abundant training data} ( $|\mathcal{S}_l^{\text{train}}|=50$). 

The proposed DRCC formulation was tested for two different radius of the Wasserstein ball.  The results of our DRCC strategy are labeled as \textit{DRCC: $\delta$}, where $\delta$ corresponds to the value of the normalized radius. Moreover, two benchmark joint optimization models were implemented denoted as: (i) \textit{SAA}, which corresponds to the proposed DRCC formulation with $\delta=0$ and thus coincides with a stochastic programming formulation solved via Sample Average Approximation, and (ii) \textit{Rob.}, which corresponds to a robust approach that solves the joint optimization problem considering the worst-case scenario for each failure time. \edit{The performance of the different models was evaluated by computing three key performance indicators (KPIs) from the outcomes collected during the simulation period:
\begin{itemize}[leftmargin=*]
    \item Percentage of preventive maintenance (\% PM): $\frac{\text{Total number of preventive maintenance activities}  }{\text{Total number of maintenance activities}}\times 100 \%$. This indicator represents the observed percentage of preventive maintenance and is expected to be high when implementing a reliable maintenance strategy. 
    \item Total operational cost: $\frac{\text{Total cumulative cost}  }{\text{Length of the simulation period}}$. This indicator represents the total operational cost per day. An effective maintenance and spare provisioning model should aim to reduce this quantity.
    \item Average chance constraint violations:
    $\frac{1}{N_t} \sum_{n=1}^{N_t} \bigl[ 0.5 \cdot \mathbb{I}(\exists$ Event violating condition of   $\mathcal{Z}_1$  in time window  $n$ ) +  $0.5 \cdot \mathbb{I}(\exists$ Event violating condition of   $\mathcal{Z}_2$  in time window  $n) \bigr]$. This indicator estimates the probability of chance constraint violations by counting and averaging the observed events that violate the conditions defining $\mathcal{Z}_1$ and $\mathcal{Z}_2$. To achieve this, it analyzes $N_t$ different realizations of the planning horizon by collecting data across $N_t$ disjoint time windows of length $T_{max}$. An effective maintenance strategy should ensure that this violation probability does not exceed the thresholds imposed by the managers.

\end{itemize}
}


\edit{The results obtained for these three KPIs are displayed in Figure \ref{fig:performance_twocases}. We first analyze the results obtained with sparse training data, indicated in red. Figure \ref{fig:performance_twocases}(b) shows that the SAA model yields the lowest operational costs, making it an effective alternative from an economic perspective. However, this approach offers poor reliability performance, resulting in solutions that may not be safe for real-world applications. Specifically, the resulting percentage of preventive maintenance (\%PM) is lower than 94\% (Figure \ref{fig:performance_twocases}(a)), and the average chance constraint violation exceeds the desired level (Figure \ref{fig:performance_twocases}(c)). This deficient reliability performance is attributed to the fact that the prognostic algorithm was trained with sparse data. Consequently, the estimated RLDs may not accurately characterize the actual RLDs, which can mislead maintenance decisions. The Rob. model, on the other hand, shows outstanding performance in terms of reliability, achieving a percentage of preventive maintenance close to 100\%. However, the Rob. model adopts a highly conservative approach to achieve this goal, which results in significantly high total operational costs.

\begin{figure}[H]
    \centering
    \includegraphics[width=0.8\textwidth]{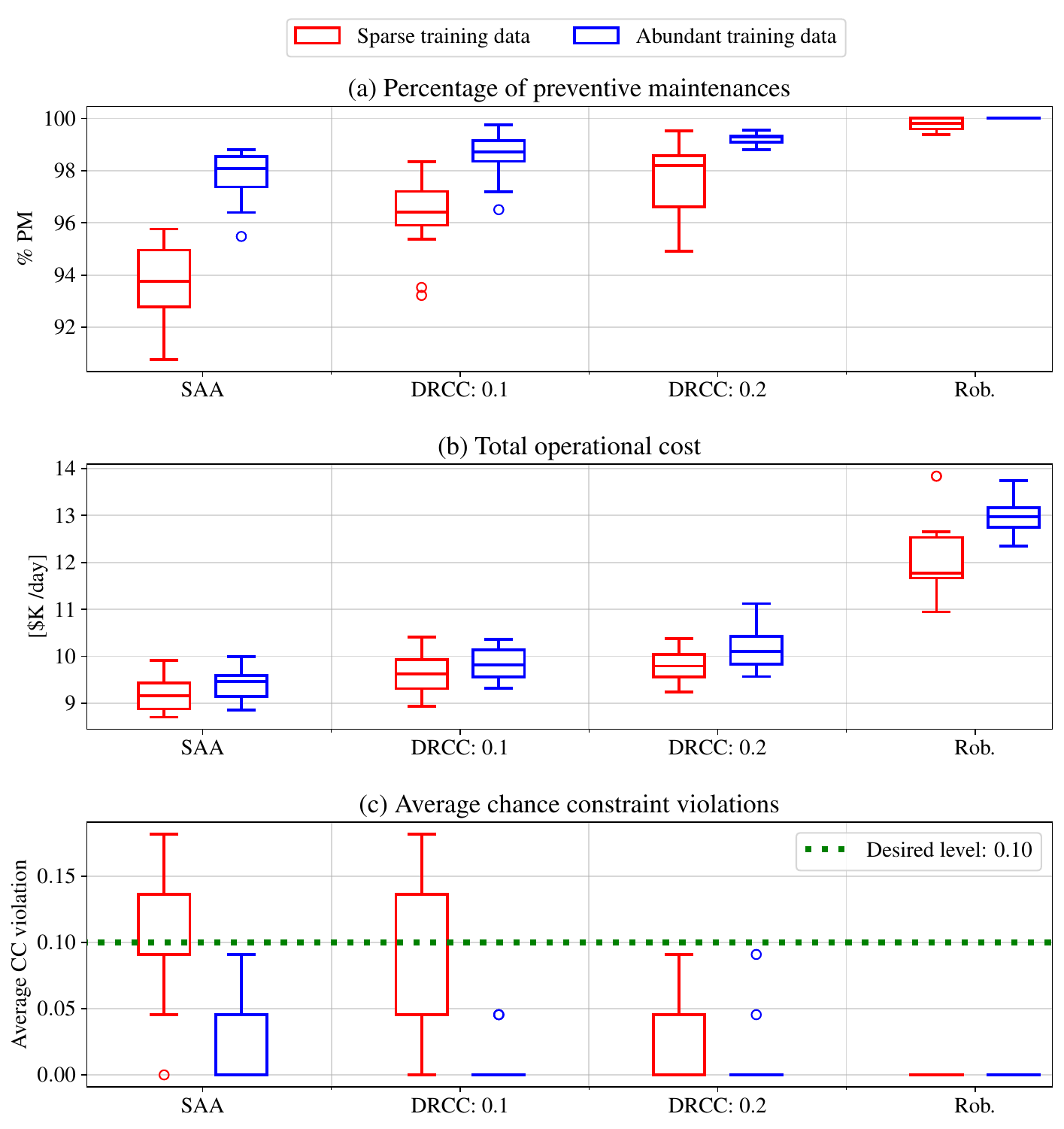}
    \caption{Results obtained considering different amounts of training data for the prognostic algorithm}
    \label{fig:performance_twocases}
\end{figure}

If the manager's primary objective is to minimize costs, the SAA model appears to be the best choice. However, our problem setting considers a scenario where system operation must comply with reliability requirements enforced by the two chance constraints. These requirements are often stipulated in operational contracts or recommended by stakeholders to enhance operational safety. When implemented with sparse training data, the SAA model may fail to meet these reliability standards, as seen in the results previously discussed. This limitation arises because the SAA model can be overconfident in its uncertainty characterization. In situations with limited training data, adopting a more cautious approach to managing uncertainty may be prudent, making robust optimization frameworks a valuable alternative to explore.
}

The performance of the proposed DRCC formulation varies depending on the value of $\delta$. As discussed in Section \ref{subsec:ambiguityset}, $\delta$ can be used to model our confidence level in the estimated RLDs. When training data is sparse, our confidence level in the prognostic results should be low. So we anticipate that increasing $\delta$ can help in seeking more robust maintenance decisions. \edit{This claim is supported by the results shown in Figure \ref{fig:performance_twocases}, which illustrates the performance of the proposed DRCC formulation for two different radii, $\delta = 0.1$ and $\delta = 0.2$, denoted as DRCC: 0.1 and DRCC: 0.2, respectively. The results indicate that increasing $\delta$ from 0.1 to 0.2 positively impacts reliability indicators. This is expected because the formulation becomes more conservative with larger $\delta$. For DRCC: 0.1, there is an average increase of 2.3\% in the percentage of preventive maintenance compared to the SAA model; however, the resulting average chance constraint violation still exceeds the desired level. Therefore, a larger $\delta$ should be considered to achieve a more conservative solution capable of satisfying the chance constraint requirements. For DRCC: 0.2, there is an average increase of 3.8\% in the percentage of preventive maintenance compared to the SAA model. More importantly, DRCC: 0.2 satisfies the desired level of chance constraints, as shown in Figure \ref{fig:performance_twocases}(c).} These improvements in terms of reliability indicators are accompanied by an increment of 0.5[\$K/day] in the total operational cost. In contrast to Rob.  model, the proposed DRCC model can satisfy chance constraint requirements without recurring to highly conservative expensive maintenance policies. As a result, our proposed DRCC model provides a reasonable balance between conservative and profitable decisions.

For the case of abundant training data (in blue), we can notice in Figure \ref{fig:performance_twocases} that the SAA model offers very satisfactory results. It shows a high percentage of preventive maintenance and is able to satisfy chance constraint requirements. This good performance is attributed to the fact that the prognostic model was trained with abundant training data and thus is able to provide accurate estimates of components' RLDs. This makes the stochastic programming model a valid alternative in this setup and gives no incentive to explore larger values for $\delta$.

\edit{\begin{rem}
 When deploying the proposed model in a rolling horizon fashion, maintenance managers have to define two application-specific parameters: $T_{max}$ and $\Delta^{upd}$. As a guideline, the selection of the planning horizon length, $T_{max}$, should consider two key aspects. 1) The planning horizon length must be shorter than the minimum guaranteed lifetime of the components under study. This ensures that Assumptions \ref{as:one_failure} and \ref{as:one_repair} are satisfied. 2) The planning horizon should be long enough to be able to evaluate the impact of placing regular spare parts at different times. On the other hand, the selection of the freeze period, $\Delta^{upd}$, depends on the system's flexibility to update decisions. Updating decisions more frequently usually results in greater efficiency \citep{rozas2024data}. However, real-world systems often require maintenance plans and team assignments to be made in advance, meaning maintenance activities cannot be adjusted too frequently. 
\end{rem}}

\begin{rem}  In practical applications, the Wasserstein radius $\delta$ can be tuned through cross-validation \citep{mohajerin2018data,xie2021distributionally}, relying on historical degradation sensor data. Specifically, cross-validation is used to evaluate the performance of the DRO model with various values of $\delta$ and then select the value that yields the best performance.
\end{rem}

\subsection{Comparison between Sequential and  Joint Optimization Models}\label{subsec:sequential_joint}

We compared the total operational costs of a sequential optimization model versus a joint optimization model for maintenance schedule and spare provisioning. To this end, we used the same simulation setup utilized for testing the case of \textit{abundant training data} in the previous section. The joint optimization model corresponds to our proposed DRCC with $\delta=0$. The sequential optimization model, on the other hand, first solves the proposed DRCC with $\delta=0$ optimizing solely maintenance schedule, i.e., omitting decision variables concerning spare provisioning. Then, it treats the resulting maintenance schedule as a known demand for spare parts to formulate the spare provisioning optimization problem.  Thus, it optimizes maintenance schedule and spare parts provisioning sequentially.

The results obtained from the simulation setup are displayed in Figure \ref{fig:joint_sequential}. On average, the sequential optimization presents total operational costs a 7\% higher than the joint optimization model. \edit{This performance difference arises because the sequential model determines the maintenance schedule without accounting for spare part inventory dynamics---its decisions are driven purely by maintenance aspects. The resulting maintenance schedule creates a deterministic demand for spare parts, which must be met by the inventory problem. This lack of coordination reduces the inventory problem's flexibility in managing spare parts orders, sometimes requiring costly expedited orders to ensure spare parts are available at the scheduled maintenance times. In contrast, the joint optimization model seamlessly integrates maintenance and spare provisioning decisions. It has the flexibility to delay or expedite repairs and spare part orders, enabling the development of cost-efficient and reliable maintenance and spare part provisioning plans. The significant difference in operational costs between the joint and sequential optimization models observed in our study highlights the importance of implementing joint optimization models for maintenance and spare provisioning.}

\begin{figure}[H]
    \centering
    \includegraphics[width=0.5\textwidth]{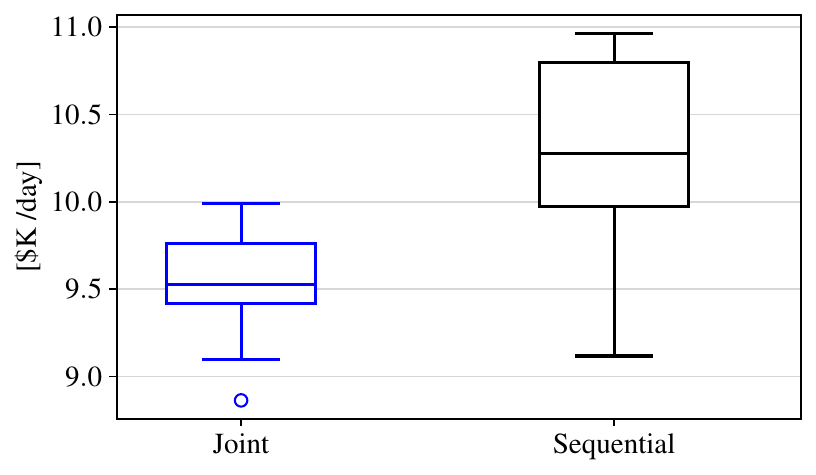}
    \caption{Comparing the total operational cost between Sequential and Joint optimization models.}
    \label{fig:joint_sequential}
\end{figure}

\subsection{Running Times}\label{subsec:running_times}

\edit{We conducted two studies to analyze the running times of the proposed DRCC joint optimization method. The first study evaluates the running time of two different solution methods for solving the proposed DRCC formulation with a fixed instance size of 150 components (30 wind turbines). The second study examines the running time of the exact solution method for four different instance sizes.}

In the first study, we evaluated the running time of two different solution methods, referred to as \textit{exact} and \textit{sample-based}, utilized to solve the proposed  DRCC formulation. The exact method fully relies on the exact reformulation derived in Section \ref{sec:formulation}. The sample-based method uses partially the exact reformulation of Section \ref{sec:formulation}, excepting the MILP reformulation of $\mathcal{Z}_2$ derived in Preposition \ref{prop:refor_cc2}. Instead, it implements the representation of $\mathcal{Z}_2$ derived in Preposition \ref{prop:cc2} via SAA according to \cite{luedtke2008sample}. The sample-based method uses $S$ independent samples of $\Tilde{\eta}$ drawn from the Bernoulli distributions specified in Proposition \ref{prop:cc2}. Then, it reformulates $\mathcal{Z}_2$ as an integer programming set that depends on the obtained samples. When solving the optimization problem using the sample-based method, the solution obtained may not be feasible to the original chance constraint because this method solves a relaxation of the original chance-constrained set. As we have an exact representation of the original chance-constrained set, we can easily verify if the solution obtained from the sample-based method denoted as $\hat{x}$ is feasible or infeasible to the original chance-constrained set. It suffices to verify if there exists some solution $(x,a,b,c)$ to \eqref{eq:final_Z2} such that $x=\hat{x}$. If we find any solution, then $\hat{x}$ is feasible to the original chance constraint; else, it is infeasible.

We solved 50 different instances of the proposed DRCC optimization, considering the two solution methods described before. \edit{Each problem instance was solved using Python 3.10 and Gurobi 9.5. The experiments were conducted on a state-of-the-art supercomputer--Partnership for an Advanced Computing Environment (PACE)}, configured with 48 GB of memory. For the sample-based method, we additionally computed the rate of chance constraint violation, which is defined as (Total chance constraint violations)/(Total number of instances). The results are presented in Table \ref{tab:running_times}. It can be seen that the sample-based method offers lower running times compared with the exact method when the number of samples is small ($S< 1000$). However, the solutions computed with the sample-based method lead to a high chance constraint violation rate, meaning that a large portion of the solutions obtained via this method are not feasible for the original problem. In order to reduce the chance constraint violation rate, the sample-based method needs to increase the number of samples. We can see that with 2000 samples, the sample-based method can solve the problem with no violations. Nevertheless, the number of samples affects the complexity of the sample-based method because this method relies on an integer programming representation that includes new binary decision variables and constraints indexed by the number of samples. We observe that when $S \geq 1000$, the running times of the sample-based method are higher than the exact method.

\edit{The proposed exact method, on the other hand, does not rely on samples, as it uses an exact reformulation of $\mathcal{Z}_2$. Table \ref{tab:running_times} shows that the exact method solves the problem with an average running time of \edit{246 seconds}. This makes it an effective tool for addressing the problem while ensuring compliance with chance constraint requirements based on the given uncertainty characterization. Consequently, we conclude that the exact method should be adopted when the goal is to make decisions while strictly enforcing the chance constraints based on the available uncertainty characterization. Conversely, if the priority were solving the problem as quickly as possible, the sample-based method might be an effective option.}

{
\renewcommand{\arraystretch}{0.7}
\begin{table}[H]
\caption{Comparing running times between the derived exact method and the sample-based approach.}
\label{tab:running_times}
\centering
\begin{tblr}{
  row{3} = {c},
  cell{1}{1} = {c},
  cell{1}{3} = {c=2}{c},
  cell{1}{5} = {r=2}{c},
  cell{2}{1} = {c},
  cell{2}{3} = {c},
  cell{2}{4} = {c},
  cell{3}{1} = {c=2}{},
  cell{4}{1} = {r=4}{c},
  cell{4}{3} = {c},
  cell{4}{4} = {c},
  cell{4}{5} = {c},
  cell{5}{3} = {c},
  cell{5}{4} = {c},
  cell{5}{5} = {c},
  cell{6}{3} = {c},
  cell{6}{4} = {c},
  cell{6}{5} = {c},
  cell{7}{3} = {c},
  cell{7}{4} = {c},
  cell{7}{5} = {c},
  hline{1} = {3-5}{},
  hline{2} = {3-4}{},
  hline{3-4,8} = {-}{},
}
      &      & Running time [ s ] &        & {Rate of chance \\constraint violation} \\
      & ~    & Average            & Max    &                                         \\
Exact &      & 246.72             & 407.38 & -                                       \\
Sample-based   & $S$=500  & 188.17             & 255.18 & 0.46                                    \\
      & $S$=1000 & 346.92             & 374.52 & 0.17                                    \\
      & $S$=1500 & 555.04             & 695.40 & 0.02                                    \\
      & $S$=2000 & 738.44             & 919.01 & 0.00                                    
\end{tblr}
\end{table}}

\edit{ The second study evaluated the running times of the exact solution method for varying instance sizes, defined by the number of wind turbines ($\mid \mathcal{K} \mid$). For each instance size, the problem was solved 50 times with different initial conditions.  Table \ref{tab:different_sizes} presents the results, including the number of decision variables and constraints for each instance size, emphasizing that the MILP reformulation derived in Section \ref{sec:formulation} leads to a large-scale MILP, even for smaller problem instances. Regarding running times, it is clear that they increase rapidly with the number of wind turbines but remain below one hour even for the largest instance of 100 turbines (300 components). In real-world applications, however, larger problem instances may need to be solved, where directly implementing the proposed model could become computationally infeasible. To address this issue, future research will focus on developing customized decomposition or approximation methods that exploit the unique structure of the problem under study.

\begin{table}[H]
\edit{
\caption{Running times of the exact method at different problem instance sizes}
\label{tab:different_sizes}
\resizebox{\textwidth}{!}{\centering
\begin{tabular}{cccccccc} 
\hline
\# Wind turbines & \# Components & \multicolumn{3}{c}{\# Decision   variables} & \multirow{2}{*}{\# Constraints} & \multicolumn{2}{c}{Running time [ s ]}  \\ 
\cline{3-5}\cline{7-8}
$\mid \mathcal{K} \mid$               & $\mid \mathcal{J} \mid$            & Continuous & Integer & Binary               &                                 & Average & Max                           \\ 
\hline
25               & 75            & 22356      & 300     & 5500                 & 93232                           & 75.86   & 98.57                         \\
50               & 150           & 59958      & 300     & 10500                & 246459                          & 246.72  & 407.38                        \\
75               & 225           & 134562     & 300     & 15500                & 545438                          & 678.93  & 1169.27                       \\
100              & 300           & 238966     & 300     & 20500                & 961817                          & 2478.68 & 3419.15                       \\
\hline
\end{tabular}}
}
\end{table}

 }
\section{Conclusions}\label{sec:conclusions}

In this paper, we have presented a DRCC formulation for jointly optimizing CBM schedules and spare parts provisioning for industrial applications. The uncertainty associated with the remaining lifetime of components was characterized by predictive analytics leveraging sensor data. Unlike previous studies, our proposed DRCC formulation acknowledges that estimated RLDs can be biased and so seeks robust solutions against distribution perturbations within a data-driven ambiguity set. This enables the joint optimization of CBM and spare parts provisioning even in presence of inaccurate prognostic results. This was showcased by the case study with wind turbines presented in Section \ref{subsec:varying_traindata}, where we could see that the hyper-parameter $\delta_j$ of our DRCC model can be adjusted in order to capture the potential perturbations of the estimated RLDs, thereby leading to more robust decisions.

We derived interesting properties of our DRCC formulation. Specifically, in Proposition \ref{prop:cc1} we proved that the DR chance constraint associated with the \edit{downtime} of critical components admits an exact integer programming representation. We provided an explicit formula to compute the parameters of that representation. In Preposition \ref{prop:cc2} we demonstrated that the DR chance constraint related to unexpected failures can be rewritten as a chance constraint whose underlying distribution is computed in a DR fashion. We also derived a closed-form expression for the parameters of that distribution. 
Combining all these results, we could show that the proposed DRCC joint optimization problem admits an exact MILP reformulation that can be efficiently solved by off-the-shelf optimizer software. In fact, in Section \ref{subsec:running_times}, we could observe that the proposed solution methodology solves the DRCC optimization problem for 150-component instances in less than 5 minutes.

Our future research will aim to extend this work including condition-based production planning. Existing literature on this topic assumes that the uncertainty on components' condition is well characterized by analytics models. However, such models can be biased due to noisy sensors or sparsity of training data. When this happens, the use of DRCC formulations may be beneficial for reducing operational costs and increasing systems' reliability. 

\section*{Acknowledgments}

This   effort  was  supported   by   NASA  under   grant  number 80NSSC19K1052 as part of the NASA Space Technology Research Institute  (STRI)  Habitats  Optimized  for  Missions  of  Exploration(HOME) ‘SmartHab’ Project. Any opinions, findings, and conclusions or recommendations expressed in this material are those of the authors and do not necessarily reflect the views of the National Aeronautics and Space Administration. Additionally, H.Rozas was supported by ANID-Chile under Fondecyt de Iniciación Grant No. 11250197 and  PIA/PUENTE AFB230002, while W. Xie was supported in part by NSF grant 2246414 and  ONR grant N00014-24-1-2066.
\addcontentsline{toc}{chapter}{Conclusion}
\setcounter{section}{0}

\section*{Appendix A. Simulation Framework}

To simulate the degradation signal of component type $l, \; l \in \mathcal{L}=\{1,2,3\}$ we utilized the following degradation model:
\begin{equation}
    S(t) = \phi \cdot \mathbb{I}(0 \leq t \leq \tau_l 
    ) + \theta_l \cdot e^{\beta_l(t-\tau_l) + \epsilon\cdot(t-\tau_l)}\cdot \mathbb{I}(t > \tau_l ),
\end{equation}
where  $\mathbb{I}(\cdot)$ is the indicator function; $\phi$ is a constant number; $\tau_l$ denotes the length of the first period and $\tau_l \sim \text{Uniform}(a_{l},b_{l})$; $\log \theta_l \sim \text{Normal}(\mu_{1,l},\sigma_{1,l}^2)$; $\log \beta_l \sim \text{Normal}(\mu_{2,l},\sigma_{2,l}^2)$, and $\epsilon$ is a centered Brownian motion such that $\epsilon(t)$ has a variance $\sigma^2_l \cdot t$. Notice that this degradation model has two phases. Phase I takes place before $\tau_l$ and is characterized by a constant signal. Phase II takes place after $\tau_l$ and models the gradual degradation process of the component. The failure event occurs when $S(t)$ crosses the failure threshold  $\Lambda_l$ for the first time.

Then degradation signals can be simulated according to the procedure described in Algorithm \ref{alg:simulation}. This procedure allows us to model the entire degradation process of components and was utilized to generate complete degradation signals. The $n$-th degradation entire degradation signal (i.e. from Phase I to failure) is denoted as $\mathcal{S}_{l}^n=\{S(0),...,S(t_{ln})\}$.

\thispagestyle{empty}
\begin{algorithm}
\caption{Simulating the Degradation Signal of component type $l$.}
\label{alg:simulation}
\begin{algorithmic}[1]
\State  $\tau_l \gets \text{Uniform}(a_{l},b_{l})$;
\State $\log \theta_l \gets \text{Normal}(\mu_{1,l},\sigma_{1,l}^2)$
\State $\log \beta_l \gets \text{Normal}(\mu_{2,l},\sigma_{2,l}^2)$
\State Initialize the time $t = 0$
\State Initialize the signal $S(t) = \phi$

\While{$S(t) \leq \Delta_l$}
    \If{$0 \leq t \leq \tau_l$}
        \State Set $S(t) = \phi$
    \Else
        \State $\epsilon \gets \text{Normal}(\mu_{l},\sigma_{l}^2)$
        \State $S(t) = \theta_l \cdot e^{\beta_l \cdot (t - \tau_l) + \epsilon\cdot(t-\tau_l)}$
    \EndIf
    \State $t=t+dt$
\EndWhile

\end{algorithmic}
\end{algorithm}
\newpage
\bibliography{references}


\end{document}